\input amstex
\documentstyle{amsppt}
\magnification=1200
\hcorrection{0.2in}
\vcorrection{-0.4in}

\define\volume{\operatorname{vol}}
\define\op#1{\operatorname{#1}}
\define\svolball#1#2{{\volume(\underline B_{#2}^{#1})}}
\define\svolann#1#2{{\volume(\underline A_{#2}^{#1})}}
\define\sball#1#2{{\underline B_{#2}^{#1}}}
\define\svolsp#1#2{{\volume(\partial \underline B_{#2}^{#1})}}

\NoRunningHeads
\NoBlackBoxes
\topmatter
\title Gromov-Hausdorff limits of Aspherical Manifolds
\endtitle
\author Xiaochun Rong \footnote{Partially supported by NSFC 11821101, BNSF Z190003, and a research fund from Capital Normal University. \hfill{$\,$}\\ 2010 Mathematics Subject Classification. Primary 53C21, 53C23, 53C24.\hfill{$\,$}}\endauthor
\address Mathematics Department, Rutgers University
New Brunswick, NJ 08903 U.S.A
\endaddress
\email rong\@math.rutgers.edu
\endemail
\address School of Mathematics Sciences, Capital Normal University, Beijing, 100048, P.R.C.
\endaddress
\abstract Let $X$ be a compact Gromov-Hausdorff limit space of a collapsing sequence of compact $n$-manifolds, $M_i$, of Ricci curvature $\op{Ric}_{M_i}\ge -(n-1)$ and all points in $M_i$ are $(\delta,\rho)$-local rewinding Reifenberg points, or sectional curvature $\op{sec}_{M_i}\ge -1$, respectively. We conjecture that if $M_i$ is an aspherical manifold of fundamental group satisfying a certain condition (e.g., a nilpotent group), then $X$ is a differentiable, or topological aspherical manifold, respectively. A main result in this paper asserts that if $M_i$ a diffeomorphic or homeomorphic to a nilmanifold, then $X$ is diffeomorphic or homeomorphic to a nilmanifold, respectively.
\endabstract
\endtopmatter
\document

\head 0. Introduction
\endhead

\vskip4mm

Consider a sequence of compact $n$-manifolds converging in the Gromov-Haudorff distance (briefly, GH-distance), $M_i@>\op{GH}>>X$, such that
$$\op{Ric}_{M_i}\ge -(n-1),\quad 0<\op{diam}(X)=d<\infty,\quad \op{vol}(M_i)\to 0.$$
In this paper, we study the following problem:

\example{Problem 0.1} Under what (topological) constraints in $M_i$, $X$ (refers as Ricci-limit space) is a manifold, differentiable, or topological.
\endexample

Note that a Ricci limit space of a non-collapsing sequence i.e, $\op{vol}(M_i)\ge v>0$,
has been well-studied: a (optimal) local regularity condition on $M_i$ is that there is $\rho>0$ such that (any) $x_i\in M_i$ is a $(\delta,\rho)$-Reifenberg point (briefly $(\delta,\rho)$-RP) i.e., for any $0<s<\rho$, the GH-distance between $s$-balls in $M_i$ centered at $x_i$ and in $\Bbb R^n$,
$$d_{\op{GH}}(B_s(x_i),\b B_s^n(0))<s\cdot \delta,\quad 0<\delta\le \delta(n).$$
According to Cheeger-Colding (\cite{CC1}), $X$ is an $n$-manifold and for $i$ large there is a diffeomorphism, $f_i: M_i\to X$, such that the push forward distance functions bi-H\"older converges. In this paper, we also use the Perel'man stability result (\cite{Per1}): if $\op{sec}{M_i}\ge -1$, then $X$ is a topological $n$-manifold, and for $i$ large there is a homeomorphism and $\epsilon_i$-GH approximation, $f_i: M_i\to X$.

Problem 0.1 is partially stimulated by the recent result of Bru\`e-Naber-Semola below (see Conjecture 9 in \cite{Za}):

\proclaim{Theorem 0.2} {\rm (\cite{BNS})} Let $M_i@>\op{GH}>>X$ be a sequence of $n$-manifolds satisfying
$$\op{sec}_{M_i}\ge -1,\quad 0<\op{diam}(X)=d<\infty,\quad \op{vol}(M_i)\to 0.$$
If $M_i$ is homeomorphic to an $n$-torus $T^n$, then $X$ is homeomorphic to $T^m$ ($m<n)$.
\endproclaim

Note that Theorem 0.2 will be
false if one replaces ``$\op{sec}_{M_i}\ge -1$'' with
``$\op{Ric}_{M_i}\ge -(n-1)$ and $\op{vol}(M_i)\ge v>0$'' (\cite{Zhu1}), or with ``an
Alexandrov metric of curvature $\ge -1$'' (\cite{BNS}).

In seeking for a general answer to Problem 0.1, we recall the following related result by Fukaya:

\proclaim{Theorem 0.3} {\rm (\cite{Fu2})} Let a sequence of aspherical $n$-manifolds (i.e., the universal covering space $\tilde M_i$ is contractible),  $M_i@>\op{GH}>>X$, such that
$$|\op{sec}_{M_i}|\le 1,\quad 0<\op{diam}(M_i)=d<\infty,\quad \op{vol}(M_i)\to 0.$$
Then $X=X'/\Gamma$ is an aspherical orbifold i.e.,  $X'$ is an aspherical manifold on which a finite group of isometries $\Gamma$ acts.
\endproclaim

According to the collapsing theory of Cheeger-Fukaya-Gromov (\cite{CG1,2}, \cite{Fu1,3}, \cite{CFG}), for $i$ large $M_i$ admits a singular nilpotent fibration (see subsection 1.3), $f_i: M_i\to X$, $f_i$ is a continuous $\epsilon_i$-GH approximation ($\epsilon_i\to 0$) such that each $f_i$-fiber, $f_i^{-1}(x), x\in X$, is an infra-nilamnifold of positive dimension, and $M_i$ admits a $C^1$-close invariant metric of high regularities. Note that in Theorem 0.3, all fibers have the same dimension (\cite{Fu2}).

Comparing Theorems 0.2 and 0.3 in dimension $n=2$: an one-parameter family of flat metrics on $T^2$, $(\epsilon S^1)\times S^1@>\op{GH}>>S^1$ ($\epsilon\to 0$), and flat metrics $g_\epsilon$ on a Klein bottle $K^2$, $(K^2,g_\epsilon)@>\op{GH}>>[0,1]\cong S^1/\Bbb Z_2$,
we introduce the following notion.

\definition{Definition 0.4} Let $\Gamma$ be a finitely generated group. We call
$\Gamma$ type $\Cal A_s$, if $\Lambda<\Gamma$ is a normal virtually
nilpotent subgroup, then any finite subgroup of $\Gamma/\Lambda$ is virtual normal i.e,
its normalizer has a finite index in $\Gamma/\Lambda$.
\enddefinition

It is easy to see that if $M$ is a compact nilmanifold i.e., $M$ is diffeomorphic to $N/\Gamma$, $N$ is a simply connected nilpotent Lie group, $\Gamma<N$ is a discrete subgroup, then $\pi_1(M)$, is type $\Cal A_s$, so is $\pi_1(M\times Y)$ type $\Cal A_s$, where $Y$
is a compact hyperbolic manifold; while $\pi_1(M')$ may not be type $\Cal A_s$ for
an infra-nilmanifold $M'$ which is not an nilmanifold ($M'$ may admit a sequence of metrics
$g_i$ of $|\op{sec}_{g_i}|\le 1$ collapsing to some $X$, not a topological manifold; e.g., $(K^2,g_\epsilon)@>\op{GH}>>[0,1]$).

In the case that $|\op{sec}_{M_i}|\le 1$, we have the following answer to Problem 0.1.

\proclaim{Theorem A} {\rm (GH-limits of aspherical manifolds)} Let a sequence of compact $n$-manifolds, $M_i@>\op{GH}>>X$, such that
$$|\sec_{M_i}|\le 1,\quad 0<\op{diam}(X)=d<\infty,\quad \op{vol}(M_i)\to 0.$$
If $M_i$ is an aspherical manifold and $\pi_1(M_i)$ is type $\Cal A_s$, then $X$ is an aspherical manifold (with a $C^\alpha$-metric tensor).
\endproclaim

Let $\pi_i: (\tilde M_i,\tilde p_i)\to (M_i,p_i)$ denote a Riemannian universal covering. Given a singular nilpotent fibration, $f_i: M_i\to X$,  it lifts to a singular nilpotennt fibration on $\tilde M_i$, $\tilde f_i: \tilde M_i\to \hat X$, where a $\tilde f_i$-fiber is a component of $\pi_i^{-1}(f_i^{-1}(x))$, $x\in X$. It turns out that a regular $\tilde f_i$-fiber is a connected nilpotent Lie group,
and the $\tilde f_i$-fibers are orbits of an $\tilde F_i$-action on $\tilde M_i$ without fixed point such that the fundamental group $\pi_1(M_i)$ preserves the $\tilde F_i$-orbits i.e., $f_i$-fibers on $M_i$ are $\pi_i$-image of $\tilde F_i$-orbits on $\tilde M_i$ (see Remark 1.14, cf. \cite{RY}).

Observe that using a $C^1$-close invariant metric on $M_i$, $\tilde F_i$ and $\pi_1(M_i)$
generate a closed group $G_i$ of isometries on $\tilde M_i$.
By now the following result implies that $X$ in Theorem A is an aspherical manifold.

\proclaim{Theorem 0.5} Let $\tilde M$ be a contractible Riemannian $n$-manifold, and let
$G$ be a closed group of isometries on $\tilde M$ such that $\tilde M/G$ is compact.

\noindent {\rm (0.5.1)} If the identity component of $G$, $G_0$, is nilpotent, then $G_0$ contains no compact subgroup.

\noindent {\rm (0.5.2)} If in (0.5.1), any finite subgroup of $G/G_0$ is virtually normal, then $G/G_0$ has no torsion element, thus $\tilde M/G$ is an aspherical manifold.
\endproclaim

In view of the above, the singular nilpotent fibration, $f_i: M_i\to X$, in Theorem A
is crucial in proving that $X$ is a manifold. Indeed, it was conjectured that a singular nilpotent fibration exists
on a collapsed manifold of curvature bound below and a local bounded covering geometry
(\cite{Ro3}, \cite{HRW}).

For $\rho>0$, let $\pi_i: (\widetilde {B_\rho(x_i)},\tilde x_i)\to (B_\rho(x_i),x_i)$ denote the Riemannian universal covering of the $\rho$-ball at $x_i\in M_i$.

\example{Conjecture 0.6} {\rm (Local bounded covering geometry, singular nilpotent fibration)} Let a sequence of compact $n$-manifolds, $M_i@>\op{GH}>>X$, $0<\op{diam}(X)=d<\infty$ and $\op{vol}(M_i)\to 0$. For $i$ large there is a singular nilpotent fibration, $f_i: M_i\to X$, under either of the following local bounded covering geometry conditions:

\noindent (0.6.1) $\op{Ric}_{M_i}\ge -(n-1)$ and $\forall\, x_i\in M_i$ is a $(\delta,\rho)$-local rewinding Reifenberg point (briefly, $(\delta,\rho)$-LRRP) i.e., $\tilde x_i$ is a $(\delta,\rho)$-RP.

\noindent (0.6.2) $\op{sec}_{M_i}\ge -1$ and $\forall\, x_i\in M_i$ is a $(\rho,v)$-local rewinding volume point (briefly, $(\rho,v)$-LRVP) i.e., $\op{vol}(B_\rho(\tilde x_i))\ge v>0$.
\endexample

Note that Conjecture 0.6 will be false, if in (0.6.1) (resp. in (0.6.2)), one weakens ``$(\delta,\rho)$-LRRP'' to ``$(\rho,v)$-LRVP'' (resp. one weakens ``$\op{sec}_{M_i}\ge -1$''
to ``$\op{Ric}_{M_i}\ge -(n-1)$''), see Example 4.2.

Clearly, Theorem 0.5 and Conjecture 0.6 suggest the following conjectured answer to Problem 0.1
i.e., Theorem A should still hold under (0.6.1), or $\op{sec}_{M_i}\ge -1$.

\example{Conjecture 0.7} (GH-limit of aspherical manifolds) Let a sequence of compact $n$-manifolds, $M_i@>\op{GH}>>X$, $0<\op{diam}(X)=d<\infty$, $\op{vol}(M_i)\to 0$, and
$$\cases \op{Ric}_{M_i}\ge -(n-1)\\ \forall\, x_i\in M_i \, \text{ is a $(\delta,\rho)$-LRRP,}\endcases \quad\text{or}\quad \sec_{M_i}\ge -1, \, \, \text{respectively}.$$
If $M_i$ is an aspherical manifold and $\pi_1(M_i)$ is type $\Cal A_s$, then $X$ is an aspherical manifold, differentiable and whose points are $(\delta,\rho')$-RP, $0<\rho'\le \rho$, or topological,
respectively.
\endexample

\proclaim{Proposition 0.8} Conjecture 0.6 implies Conjecture 0.7.
\endproclaim

Note that ``$(\rho,v)$-LRVP'' in (0.6.2) is not assumed in Conjecture 0.7; which
was implied by the following result of Gromov.

\proclaim{Theorem 0.9} {\rm (\cite{Gr2}, cf. \cite{BNS})} Let $M$ be a compact aspherical $n$-manifold of $\op{Ric}_M\ge -(n-1)$ and $\op{diam}(M)\le d$. Then for $\tilde x\in \tilde M$, $\op{vol}(B_1(\tilde x))\ge v(n,d)>0$, a constant depending on $n$ and $d$.
\endproclaim

We now state the main result in this paper.

\proclaim{Theorem B} {\rm (Conjecture 0.7 in nilmanifolds)} Let a sequence of $n$-manifolds, $M_i@>\op{GH}>>X$, $0<\op{diam}(X)=d<\infty$ and $\op{vol}(M_i)\to 0$.

\noindent {\rm (B1)} If $M_i$ is a nilmanifold of $\op{Ric}_{M_i}\ge -(n-1)$ and $\forall\, x_i\in M_i$ is $(\rho,\delta)$-LRRP, then $X$ is diffeomorphic to a nilmanifold whose points are
$(\delta,\rho')$-RP, $0<\rho'\le \rho$.

\noindent {\rm (B2)} If $M_i$ of $\op{sec}_{M_i}\ge -1$ is homeomorphic to a nilmanifold, then $X$ is homeomorphic to a nilmanifold.
\endproclaim

Note that (B1) will be false if one replaces (local bounded covering geometry) ``$x_i\in M_i$ is a $(\delta,\rho)$-LRRP'' by ``$x_i$ is a $(\rho,v)$-LRVP'' (\cite{Zhu1}), (B2) generalizes Theorem 0.2 and also gives an affirmative answer to a problem of A. Petrunin \footnote{In a
private communication with Zamora about Theorem B, the author was told that Petrunin
asked (B2) as a problem concerning \cite{BNS}.}, and (B2) will be false if one replaces ``$\op{sec}_{M_i}\ge -1$'' by an Alexandrov metric of curvature $\ge -1$ (\cite{BNS}).

Because Conjecture 0.6 has been wild open for $X$ a singular space, let's consider the following commutative diagram on an equivariant GH-convergence (passing to a subsequence):
$$\CD (\tilde M_i,\tilde p_i,\Gamma_i)@>\op{eqGH}>>(\tilde X,\tilde p,G)\\ @VV \pi_i V @VV /G V\\
(M_i,p_i)@>\op{GH}>>(X=\tilde X/G,p),\endCD\tag 0.10$$
where $\Gamma_i=\pi_1(M_i,p_i)$, and $G$ is the limit group of isometries on $\tilde X$ (a Lie group, \cite{CoN}). Obviously, $X$ is a manifold, if $\tilde X$ is an $n$-manifold and the $G$-action is free.

By Theorem 0.9 and stability results mentioned earlier (\cite {CC1}, \cite{Per1}), one concludes that in Conjecture 0.7, $\tilde X$ (in (0.10)) is an $n$-manifold, differentiable or topological, respectively. In the proof of Theorem B, a main technical result is the following.

\proclaim{Theorem 0.11} {\rm (Contractibility)} Let the assumptions be as in Theorem B. Then $\tilde X$ in (0.10)
is a contractible $n$-manifold, differentiable or topological, respectively.
\endproclaim

A consequence of Theorem 0.11 and a version of Theorem 0.5 for GH-limit spaces (see Lemma 3.5) is that $X$ is an asphercial manifold, differential or topological, respectively. Because $G$ is nilpotent, $\pi_1(X)$ is a nilpotent group, thus $X$ is homotopy equivalent to a nilmanifold. By now (B2) follows from a result of Farrell-Hsiang (\cite{FH}): a homotopy nilmanifold is homeomorphic to a nilmanifold (see Theorem 3.6).

In (B1), to show that $X$ is diffeomorphic to a nilmanifold, the following partial verification of (0.6.1) is unexpendable.

\proclaim{Theorem 0.12} {\rm (Nilpotent fibration, \cite{Ro4}, see Remark 1.17)} Let $M_i@>\op{GH}>>X$ be a collapsing sequence of compact $n$-manifolds satisfying
$$\op{Ric}_{M_i}\ge -(n-1),\quad 0<\op{diam}(X)=d<\infty,\quad \forall\, x_i\in M_i \text{ is a $(\rho,v)$-LRVP}.$$
If all points in $X$ are $(\delta,\rho)$-RP, then for $i$ large there is a smooth fiber bundle, $F_i\to M_i@>f_i>>X$, with $F_i$ an infra-nilmanifold and affine structural group, and $f_i$ is an $\epsilon_i$-Gromov-Hausdorff approximation, $\epsilon_i\to 0$.
\endproclaim

First, it is easy to check that in (B1), points in $\tilde X$ are $(\delta,\tilde \rho)$-RP ($0<\tilde \rho<\rho$), and points in $\tilde X/G$ are $(\delta,\rho')$-RP (see Proposition 1.11). Applying Theorem 0.12, for $i$ large one
gets a fiber bundle, $F_i\to M_i@>f_i>>X$, with $F_i$ a nilmanifold (because $M_i$ is nilmanifold) and affine structural group. Consequently, the lifting fiber bundle on $\tilde M_i$ is a principal $\tilde F_i$-bundle
$\tilde F_i\to \tilde M_i@>\tilde f_i>>\hat X$, with $\tilde F_i$ a simply connected nilpotent
Lie group and $\hat X$ is universal covering space of $X$. Because $X$ is aspherical,
$\hat X$ is contractible, the principal $\tilde F_i$-bundle is a trivial bundle.
Because the inclusion, $F_i\hookrightarrow M_i$, induces injection $\pi_1(F_i)\hookrightarrow \pi_1(M_i)$, by Malc\'ev rigidity we identify $\tilde F_i$ with a normal subgroup of
$\tilde M_i$, thus $X$ is diffeomorphic to the quotient of simply connected nilpotent Lie group, $\tilde M_i/\tilde F_i$, by a co-compact discrete subgroup, $\pi_1(M_i)/\pi_1(F_i)$.

Our proof of Theorem 0.11 relies on the property that a simply connected nilpotent Lie group $N$ is exponential i.e., the exponential map from the Lie algebra, $\op{Exp}: T_eN\to N$, is a diffeomorphism (see Theorem 3.1).

Nevertheless, the following is possible (see Remark 3.4).

\example{Conjecture 0.13} (Contractibility) Let a sequence of compact aspherical $n$-manifolds, $M_i@>\op{GH}>>X$,
$$\cases \op{Ric}_{M_i}\ge -(n-1)\\ \forall\, x_i\in M_i \, \text{ is a $(\delta,\rho)$-LRRP,}\endcases \quad\text{or}\quad \sec_{M_i}\ge -1, \, \, \text{respectively}.$$
Then $\tilde X$ in (0.10) is a contractible $n$-manifold, differentiable and
points are $(\delta,\rho')$-RP, or topological, respectively.
\endexample

\proclaim{Proposition 0.14} Conjecture 0.13 implies Conjecture 0.7.
\endproclaim

In dimension $3$, we have the following result which implies Conjecture 0.7 in dimension $3$.

\proclaim{Theorem C} Let a sequence of asphercial $3$-manifolds, $M_i@>\op{GH}>>X$, such that $$\op{Ric}_{M_i}\ge -2,\quad 0<\op{diam}(X)=d<\infty, \quad \op{vol}(M_i)\to 0.$$
If $\pi_1(M_i)$ is type $\Cal A_s$, then $X$ is an aspherical manifold.
\endproclaim

Theorem C generalizes a result in \cite{BNS}: Theorem 0.2 holds in dimension $3$ under $\op{Ric}_{M_i}\ge -2$; Theorem C applies to $M\in \{S^1\times \Sigma_g, g\ge 1\}$, where $\Sigma_g$ denotes a closed surface of genus $g$ in $\Bbb R^3$, and $(M,g_\epsilon)=(\epsilon S^1)\times \Sigma_g@>\op{GH}>>\Sigma_g$ ($\epsilon\to 0$) with bounded sectional curvature.

\vskip4mm

This rest paper is organized as follows:

In Section 1, we will review notions and results that will be used throughout this paper.

In Section 2, we will prove Theorem 0.5, Theorem A, and Proposition 0.8

In Section 3, we will prove Theorems 0.11, Theorem B, and Proposition 0.13.

In Section 4, we will prove Theorem C.

\vskip2mm

{\bf Acknowledgement}: The author would like to thank Jiayin Pan, Jikang Wang, and Sergio Zamora for useful comments on this preprint, and thank an anonymous referee for bring up \cite{Fu2} to the author's attention.

\vskip4mm

\head 1. Preliminaries
\endhead

\vskip4mm

For convenience of readers, we will give a brief review of notions and results that are used in
the rest of paper.

\vskip4mm

\subhead 1.1. Equivariant Gromov-Hausdorff convergence
\endsubhead

\vskip4mm

Let $X, Y$ be compact metric spaces. For $\epsilon>0$, a map, $f: X\to Y$, is called an
$\epsilon$-Gromov-Hausdorff approximation (briefly, GHA), if $f$ is an $\epsilon$-isometry i.e., $|d_Y(f(x),f(x'))-d_X(x,x')|<\epsilon$ and $f(X)\subseteq Y$ is $\epsilon$-dense.

We say that a sequence of compact metric spaces, $X_i$, Gromov-Hausdorff (briefly, GH) converges to a compact metric space, denoted by $X_i@>\op{GH}>>X$, if there is a sequence of $\epsilon_i$-GHA, $f_i: X_i\to X$, $\epsilon_i\to 0$.

A sequence of compact metric spaces, $X_i$, contains a GH-convergent subsequence if the following Ascoli type precompactness conditions are satisfied:

\noindent (1.1.1) (Uniform bound) the diameter, $\op{diam}(X_i)\le d$,

\noindent (1.1.2) (Equi-continuous) For $\epsilon>0$, the number of points in an $\epsilon$-net of $X_i$, $|X_i(\epsilon)|\le c(\epsilon)$, a constant depending on $\epsilon$ (an $\epsilon$-net of $X$ is a subset of points that are $\epsilon$-separated and $\epsilon$-dense).

There is a natural metric on the isometry group of $X$, for $\alpha, \beta\in \op{Isom}(X)$,
$d(\alpha,\beta)=\max\{d(\alpha(x),\beta(x)),\, x\in X\}$, and it is easy to check
that $(\op{Isom}(X),d)$ is a compact metric space.

For a sequence of pairs, $(X_i,\Gamma_i)$, and a pair $(X,\Gamma)$, where $\Gamma_i, \Gamma$ are
closed subgroups of $\op{Isom}(X_i)$, we say that $(X_i,\Gamma_i)$ equivariant GH-convergent to
$(X,\Gamma)$, denoted by $(X_i,\Gamma_i)@>\op{eqGH}>>(X,\Gamma)$, if there is a sequence of pairs, $(f_i, \phi_i): (X_i,\Gamma_i)\to (X,\Gamma)$, $f_i: X_i\to X$ and $\phi_i: \Gamma_i\to \Gamma$ are $\epsilon_i$-GHA's, such that
$$d_X(f_i(\gamma_i(x_i)), \phi_i(\gamma_i)(f_i(x_i))<\epsilon_i,\quad x_i\in X_i, \gamma_i\in \Gamma_i.$$

\proclaim{Lemma 1.2} Assume that a sequence of compact metric spaces, $X_i@>\op{GH}>>X$. Let $\Gamma_i<\op{Isom}(X_i)$ be a closed subgroup. Then passing to a subsequence, $(X_i,\Gamma_i)@>\op{eqGH}>>(X,\Gamma)$,
where $\Gamma$ is a closed subgroup of $\op{Isom}(X)$. Moreover, the sequence of compact orbit spaces equipped with quotient metrics, $X_i/\Gamma_i@>\op{GH}>>X/\Gamma$.
\endproclaim

Next, we review notions that extend the above GH-convergence and equivariant GH-convergence to a sequence of complete locally compact metric spaces (which may not be compact).

For $\epsilon>0$, a map between two pointed metric spaces, $f: (X,p)\to (Y,q)$, $q=f(p)$, is called an $\epsilon$-GHA, if $f: B_{\epsilon^{-1}}(p)\to B_{\epsilon^{-1}+\epsilon}(q)$ is an $\epsilon$-GHA.

A sequence of pointed complete locally compact metric spaces, $(X_i,p_i)$, GH-converges to a pointed metric space $(X,p)$, denoted by $(X_i,p_i)@>\op{GH}>>(X,p)$, if there is a sequence of $\epsilon_i$-GHA, $f_i: (X_i,p_i)\to (X,p)$, $f_i(p_i)=p$, $\epsilon_i\to 0$. If $\Gamma_i$ is a closed subgroup of $\op{Isom}(X_i)$, and $\Gamma$ a closed subgroup of $\op{Isom}(X)$, a pair of maps, $(f_i,\phi_i): (X_i,p_i,\Gamma_i)\to (X,p,\Gamma)$, is called an $\epsilon_i$-eqGHA, if $f_i: (X_i,p_i)\to (X,p)$ is an $\epsilon_i$-GHA, and $\phi_i: \Gamma_i(\epsilon_i)\to \Gamma(\epsilon_i)$, such that
$$d_X(f_i(\gamma_i(x_i), \phi_i(\gamma_i)(f_i(x_i))<\epsilon_i,\quad x_i\in B_{\epsilon_i^{-1}}(p_i), \gamma_i\in \Gamma_i, \gamma_i(x_i)\in B_{\epsilon_i^{-1}}(p_i),$$
where $\Gamma_i(\epsilon)=\{\alpha_i\in \Gamma_i, \, d_i(\tilde p_i,\alpha_i(p_i))<\frac 12\epsilon^{-1}\}$.

A metric space $X$ is called a length space, if for $p, q\in X$, $d(p,q)$ is realized by
the length of a path from $p$ to $q$.

\proclaim{Lemma 1.3} Assume that $(X_i,p_i)$ are complete locally compact length metric spaces.

\noindent {\rm (1.3.1)} That $(X_i,p_i)@>\op{GH}>>(X,p)$ if and only if for any $r>0$, the closed balls, $(\bar B_r(p_i),p_i)@>\op{GH}>>(\bar B_r(p),p)$.

\noindent {\rm (1.3.2)} Assume that $(X_i,p_i)@>\op{GH}>>(X,p)$, and closed subgroup $\Gamma_i<\op{Isom}(X_i)$. Then passing to a subsequence, $(X_i, p_i,\Gamma_i)@>\op{eqGH}>>(X,p,\Gamma)$, where $\Gamma$ is a closed subgroup of $\op{Isom}(X)$.
Moreover, $(X_i/\Gamma_i,\bar p_i)@>\op{GH}>>(X/\Gamma,\bar p)$.

\noindent {\rm (1.3.3)} (see Appendix) Let $\Lambda_i<\Gamma_i$ denote a normal closed subgroup of $\Gamma_i$, such that $\Lambda_i\to G_1$. Then the following diagram commutes:
$$\CD (X_i,p_i,\Lambda_i<\Gamma_i)@>\op{eqGH}>>(X,p,G_1<G)\\
@VV /\Lambda_i V @VV /G_1V\\
(X/\Lambda_i,\bar p_i,\Gamma_i/\Lambda_i)@>\op{eqGH}>>(X/G_1, \bar p,G/G_1)\endCD$$
\endproclaim

For a reference of this subsection, see \cite{FY}, \cite{Ro2}.

\vskip4mm

\subhead 1.2. Cheeger-Colding-Naber theory on Ricci limit spaces
\endsubhead

\vskip4mm

Let $M$ be a complete $n$-manifold of $\op{Ric}_M\ge (n-1)H$, $H$ is a constant. A basic geometric property is the Bishop-Gromov relative volume comparison; where Gromov
extended Bishop volume comparison, aiming for obtain the following precompactness
which has been a fundamental platform in metric Riemannian geometry.

\proclaim{Theorem 1.4} {\rm(Precompactness, Gromov)} Any sequence of complete $n$-manifolds $M_i$ of $\op{Ric}_{M_i}\ge -(n-1)$ contains a GH-convergent subsequence, $(M_i,p_i)@>\op{GH}>>(X,p)$.
\endproclaim

In literature, $X$ in Theorem 2.1 is referred as a Ricci limit space, which is a length metric space. An application of Theorem 1.4 is that for $x\in X$, and any sequence of reals, $r_i\to \infty$, there is a subsequence such that $(r_{i_k}X,x)@>\op{GH}>>(C_xX,x)$, a tangent cone at $x$, which may depend on a choice of a subsequence.

Since 1990's, there has been intensive study in interplays between geometric and topological structures on $M_i$ and on $X$, and the Cheeger-Colding-Naber theory on Ricci limit spaces provides basic tools (\cite{CC1-3}, \cite{Co1,2}, \cite {ChN}, \cite{CoN}, \cite{CJN}, etc), and a fairly good understanding has been obtained for a non-collapsing sequence i.e., $\op{vol}(B_1(p_i))\ge v>0$.  For instance, $X$ has an open dense subset that is an $n$-manifold, whose complement has Hausdorff co-dimension at least $2$, every tangent
cone is an $n$-dimensional metric cone, etc.

For a non-collapsing sequence, a fundamental  result is the following stability theorem of Cheeger-Colding.

\proclaim{Theorem 1.5} {\rm (Stability, \cite{CC1}, \cite{CJN})} Let a sequence of complete $n$-manifolds, $(M_i,p_i)@>\op{GH}>>(X,p)$, such that $\op{Ric}_{M_i}\ge -(n-1)$ and $\forall \, x_i\in M_i$ is $(\delta,\rho)$-RP, $0<\delta\le \delta(n)$.

\noindent {\rm (1.5.1)} If $0<\op{diam}(X)=d<\infty$, then $X$ is an $n$-manifold, and for $i$ large, there is a diffeomorphism, $f_i:  X\to M_i$, such that the pullback distance functions bi-H\"older converges.

\noindent {\rm (1.5.2)} If $X$ is not compact, then fixing $r>0$, for $i$ large there is a diffeomorphism, $f_i: B_r(p_i)\to B_r(p)$ (open balls), such that the pullback distance functions bi-H\"older converges converges.
\endproclaim

The following stability result of Perelman will also be used in this paper; for simplicity here we will state it for Riemannian manifolds.

\proclaim{Theorem 1.6} {\rm (Stability, \cite{Per1}, \cite{Ka2})} Let a sequence of complete $n$-manifolds, $(M_i,p_i)@>\op{GH}>>(X,p)$, such that $\op{sec}_{M_i}\ge -1$ and $\op{vol}(B_1(p_i))\ge v>0$.

\noindent {\rm (1.6.1)} If $0<\op{diam}(X)=d<\infty$, then $X$ is a topology $n$-manifold, and for $i$ large, there is a homeomorphism, $f_i: M_i\to X$, which is also $\epsilon_i$-GHA, $\epsilon_i\to 0$.

\noindent {\rm (1.6.2)} If $X$ is not compact, then fixing $r>0$, for $i$ large there is a homeomorphism, $f_i: B_r(p_i)\to B_r(p)$, which is also $\epsilon_i$-GHA, $\epsilon_i\to 0$.
\endproclaim

In contrast with our understanding for a noncollapsing sequence, very limited knowledge has been gained for a collapsing sequence (see Theorems 1.7-1.9 below). On the other hand, recent examples have revealed
properties of $X$ as `bad' as one could imagine; the Hausdorff dimension $\dim_H(X)$ may not be an integer (\cite{PW}, \cite{Pan}), and $X$ may not contain a single topological manifold point (\cite{HNW}, \cite{Zh2}).

\proclaim{Theorem 1.7} {\rm (\cite{CoN})} The isometry group of a Ricci limit space $X$ is a Lie group. In particular, any closed group of isometries on $X$ is a Lie group.
\endproclaim

Theorem 1.7 is required in Kapovitch-Wilking's resolution on the following conjecture of Gromov.

\proclaim{Theorem 1.8} {\rm (Generalized Margulis lemma, \cite{KW})} For $n\ge 2$, there are
constants $\epsilon(n), c(n)>0$, such that for any complete $n$-manifold $M$ of $\op{Ric}_M\ge -(n-1)$, and $x\in M$, $\Lambda_x=\op{Im}[\pi_1(B_{\epsilon(n)}(x))\to \pi_1(B_1(x))]$ is $c(n)$-nilpotent i.e., $\Lambda_x$ contains a nilpotent subgroup of index $\le c(n)$.
\endproclaim

Theorem 1.7 is also used in Wang's result on a topological property of a Ricci limit space.

\proclaim{Theorem 1.9} {\rm (Semi-locally simply connectedness, \cite{Wa})} A Ricci limit space $X$ is semi-locally simply connected. In particular, the universal covering space of $X$ is simply connected.
\endproclaim

We now state a list of properties which will be used in the rest of the paper, and which can be easily derived from the notions and results in the above.

\proclaim{Proposition 1.10} Let a sequence of compact $n$-manifolds of $\op{Ric}_{M_i}\ge -(n-1)$, $M_i@>\op{GH}>>X$ ($X$ is compact), and let $\pi_i: (\tilde M_i,\tilde p_i)\to (M_i,p_i)$ denote the Riemannian universal covering space. Then passing to a subsequence, the following equivariant GH-convergence diagram commutes:
$$\CD (\tilde M_i,\tilde p_i,\Gamma_i)@>\op{eqGH}>>(\tilde X,\tilde p,G)\\
@VV \pi_i V @VV /GV \\
(M_i,p_i)@>\op{GH}>>(X=\tilde X/G,p),\endCD$$
where $\Gamma_i=\pi_1(M_i,p_i)$.

\noindent {\rm (1.10.1)} $G$ is a Lie group, and the identity component $G_0$ is nilpotent.

\noindent {\rm (1.10.2)} There is $\epsilon>0$, such that the subgroup, $\Gamma_i(\epsilon)$,
is normal and contains a nilpotent subgroup of index $\le c(n)$, $\Gamma_i(\epsilon)\to G_0$, and for $i$ large, $\Gamma_i/\Gamma_i(\epsilon)\cong G/G_0$, where $\Gamma_i(\epsilon)$ is the subgroup generated by elements of $\Gamma_i$ whose maximal displacement in $B_1(\tilde p_i)$ is less than $\epsilon$.

\noindent {\rm (1.10.4)} The universal covering of $X$ is $(\tilde X/G_0)/H$, where $H<G/G_0$ is a subgroup generated by isotropy groups of the $G/G_0$-action on $\tilde X/G_0$
\endproclaim

Finally, we supply a result mentioned in the introduction following Theorem 0.12.

\proclaim{Proposition 1.11} Let a sequence of compact $n$-manifolds, $M_i@>\op{GH}>>X$, satisfying
$$\CD (\tilde M_i,\tilde p_i,\Gamma_i)@>\op{eqGH}>>(\tilde X,\tilde p,G)\\
@VV \pi_i V @VV /GV \\
(M_i,p_i)@>\op{GH}>>(X=\tilde X/G,p),\endCD \quad \cases \op{Ric}_{M_i}\ge -(n-1) \\ 0<\op{diam}(X)=d<\infty.\endcases$$

\noindent {\rm (1.11.1)} Assume that $x_i\in M_i \text{ is a $(\rho,v)$-LRVP}$ and $x\in X \text{ is a $(\delta,\rho)$-RP}$. Then $x_i\in M_i$ is a $(\delta,\rho')$-LRRP, where $0<\rho'=\rho'(n,d,v)<\rho$.

\noindent {\rm (1.11.2)} Assume that $x_i\in M_i \text{ is a $(\delta,\rho)$-LRRP}$ and $\op{vol}(B_1(\tilde p_i))\ge v>0$. Then $\tilde x_i\in \tilde M_i$ is a $(\delta,\rho')$-RP, where $0<\rho'=\rho'(n,d,v)<\rho$.
\endproclaim

\demo{Proof} (1.11.1) For $x_i\in M_i$, $x_i$ is a $(\rho,v)$-LRVP i.e.,  if $\pi_i': (\widetilde{B_\rho(x_i)},\tilde x_i')\to (B_\rho(x_i),x_i)$ denotes the Riemannian universal covering space, then $\op{vol}(B_\rho(\tilde x_i'))\ge v>0$, and if $x_i\to x\in X$, then
$x$ is $(\delta,\rho)$-RP. We shall show that for $i$ large $x_i$ is a $(\delta,\rho')$-LRRP i.e., $\tilde x_i'$ is a $(\delta,\rho')$-RP, for some $0<\rho'<\rho$.

Let $\Lambda_i=\pi_1(B_\rho(x_i),x_i)$. Using a domain version of pre-compactness theorem in \cite{Xu}, we may assume the following commutative diagram (where points in consideration are away from
$(\pi_i')^{-1}(\partial B_\rho(x_i))$):
$$\CD (\widetilde{B_\rho(x_i)},\tilde x_i',\Lambda_i)@>\op{eqGH}>>(\tilde Y,\tilde x,Q)\\
@VV \pi_i' V @VV \op{proj} V \\
(B_\rho(x_i),x_i)@>\op{GH}>>(B_\rho(x)=\tilde Y/Q\subset X,x).\endCD$$
Because $B_\rho(x)=\tilde Y/Q$ is a smooth manifold and points in $B_{\frac \rho2}(x)$ are
$(\delta,\frac \rho2)$-RP, $\tilde Y$ is a smooth manifold and the $Q$-action is free
($Q$ is a Lie group).

We shall show that for $i$ large $\tilde x_i'$ is a $(\delta,\rho')$-RP via
the following criterion (cf. \cite{HKRX}): let $M$ be an $n$-manifold of $\op{Ric}_M\ge -\delta(n-1)$, $0<\delta\le \delta(n)$, $x\in M$ is a $(\delta,\rho)$-RP if and only if $\op{vol}(B_\rho(x))\ge (1-\delta)\op{vol}(\b B_\rho(\Bbb R^n))$.
Because $\tilde x_i'\to \tilde x$, by Theorem 1.5, it suffices to show that $\op{vol}(B_{\rho'}(\tilde x))\ge (1-\delta)\op{vo}(\b B_{\rho'}(\Bbb R^n))$.

To see a desired volume estimate based on that $x$ is $(\delta,\rho)$-RP, our approach is to magnify around $\tilde x$: assuming $r_j\to \infty$ such that $(r_jX, y)\to (C_yX,y)$, and passing to a subsequence, we may assume the following commutative diagram:
$$\CD (r_jB_{\rho'}(\tilde x),\tilde x,Q_0)@>\op{eqGH}>>(C_{\tilde x}\tilde Y\cong \Bbb R^k\times \tilde Z,\tilde x,\Bbb R^k)\\
@VV \op{proj} V @VV \op{proj} V \\
(r_j B_{\rho'}(x),x)@>\op{GH}>>(C_xX=\tilde Z,x),\endCD$$
where $\rho'>0$ is taken so that any element in $Q/Q_0$ has a maximal displacement $\ge \rho'$ in $B_{\frac \rho2}(\tilde x)$, and $\Bbb R^k$ denotes the tangent space of $Q_0(\tilde x)$ at $\tilde x$ (the splitting property of Ricci limit spaces). Note that $\op{proj}: \tilde Z\to C_xX$ is an isometry. Because
$x$ is a $(\delta,\rho)$-RP, by now it is clear that for $r_j$ large, or $\rho'$ suitably small, $\op{vol}(B_{\rho'}(\tilde x))\ge (1-\delta)\op{vol}(\b B_{\rho'}(\Bbb R^n))$.

(1.11.2) Let $\pi_i: (\tilde M_i,\tilde x_i)\to (M_i,x_i)$ denote the Riemannian universal covering map. Observe that restricting to connected components, $\hat \pi_i: ((\pi_i')^{-1}(B_\rho(x_i)),\tilde x_i')\to (\pi_i^{-1}(B_\rho(x_i)),\tilde x_i)$
is a covering map. By Bishop-Gromov relative volume comparison and that $\op{vol}(B_1(\tilde p_i))\ge v>0$, we derive
$$|\hat \pi_i^{-1}(\tilde x_i)\cap B_{\frac \rho4}(\tilde x_i')|\le
\frac{\op{vol}(B_\rho(\tilde x_i'))}{\op{vol}(B_{\frac\rho4}(\tilde x_i))}\le \frac{\op{vol}(\b B_\rho(\Bbb H^n))}{c(n,\rho,v)}\le c'(n,d,\rho,v),$$
where $\Bbb H^n$ denotes the simply connected $n$-dimensional hyperbolic space.
We claim that for any $\gamma_i\in \pi_1(B_\rho(\tilde x_i))$, the orbit of the cyclic subgroup generated by $\gamma_i$, $\left<\gamma_i\right>(\tilde x_i')$, is not contained in $B_{\frac \rho4}(\tilde x_i')$. Consequently, the minimal displacement of $\gamma_i$ at $\tilde x'$ is $\ge \rho'=\rho'(n,d,\rho,v)>0$, thus $\tilde x_i$ is a $(\delta,\rho')$-RP.

If the claim is false, because $B_\rho(\tilde x_i')$ is diffeomorphic and bi-H\"older close to $\b B^n_\rho(0)$, using the Euclidean metric the energy function, $e(\cdot)=\sum_{j=1}^{o(\gamma_i)}\frac 12\b d^2_0(\gamma_i^j(\tilde x_i'),\cdot)$, achieves the minimum at unique point $\tilde z_i\in B_{\frac \rho2}(\tilde x_i')$, where the order, $o(\gamma_i)\le c'(n,d,\rho,v)$. Consequently, $\left<\gamma_i\right>(\tilde z_i)$ concentrates around $\tilde z_i$, a contradiction.
\qed\enddemo

\vskip4mm

\subhead 1.3. Collapsed manifolds of local Ricci bounded covering geometry
\endsubhead

\vskip4mm

In ths section, we will briefly describe a background of Theorem 0.12, and an outline of
proof.

An $n$-manifold $M$ with a (normalized) bound on curvature (to prevent a possible scaling)
is called $\epsilon$-collapsed, if for any $x\in M$, $\op{vol}(B_1(x))<\epsilon$. Examples
of collapsed $n$-manifolds of various curvature bounds,
$$|\op{sec}_M|\le 1, \quad \op{sec}_M\ge -1,\quad |\op{Ric}_M|\le n-1,\quad \op{Ric}_M\ge -(n-1),$$
indicate that the complexity in the underlying geometric and topological structure of $M$
consistent with the regularity strength of a curvature bound; roughly speaking, geometric and topological structures on a collapsed manifold $M$ of $|\op{sec}_M|\le 1$ is well-understood, while collapsed manifolds under other curvature bounds are less (or far) from understood.

According to the collapsing theory of Cheeger-Fukaya-Gromov (\cite{Gr1}, \cite{CG1,2}, \cite{Fu1, 3}, \cite{CFG}), an $\epsilon$-collapsed $n$-manifold $M$ of $|\op{sec}_M|\le 1$ admits a mixed nilpotent structure (\cite{CFG}) and a $C^1$-close
invariant metric of high regularities. Conversely, if a manifold $M$ admits a topological mixed nilpotent structure, then $M$ admits a one-parameter family of $\epsilon$-collapsed metrics $g_\epsilon$ of $|\op{sec}_{g_\epsilon}|\le 1$, $\epsilon\to 0$ (\cite{CR}).

In the Cheeger-Fukaya-Gromov theory, a corner stone is the following.

\proclaim{Theorem 1.12} {\rm (Almost flat manifolds, \cite{Gr1}, \cite{Ruh})}  For $n\ge 2$, there are constants, $\epsilon(n), c(n)>0$, such that if an $n$-manifold $M$ of $|\op{sec}_M|\le 1$ satisfies that $\op{diam}(M)<\epsilon(n)$, then $M$ is diffeomorphic to an infra-nilmanifold, $N/\Gamma$, where $N$ is a simply connected nilpotent Lie group, $\Gamma<N\rtimes \op{Aut}(N)$
is a discrete subgroup, $[\Gamma:\Gamma\cap N]\le c(n)$.
\endproclaim

The following can be treated as a parametrized Theorem 1.12.

\proclaim{Theorem 1.13} {\rm (Singular nilpotent fibrations, \cite{Fu1,3}, \cite{CFG})} Assume a sequence of compact $n$-manifolds converging in the Gromov-Hausdorff topology, $M_i@>\op{GH}>> X$, such that
$$|\op{sec}_{M_i}|\le 1,\quad 0<\op{diam}(X)=d, \quad \dim(X)<n.$$
Then for $i$ large, there is an $\op{O(n)}$-invariant fiber bundle on the frame bundle $F(M_i)$ equipped with a canonical metric, such that the following diagram commute,
$$\CD \hat F_i\to (F(M_i),\op{O(n)})@>\hat f>>(Y,\op{O(n))}\\
@VV\phi_i V @VV/\op{O(n)}V\\
M_i@>f_i>>X=Y/\op{O(n)},\endCD\tag 1.13.1$$
where $\hat F_i$ is a nilmanifold and the structural group is affine, and $\phi_i(\hat F_i)$ is an infra-nimanifold of a positive dimension, such that $Y$ is a $C^{1,\alpha}$-manifold and
$\hat f_i$ satisfies the following regularities:

\noindent {\rm (1.13.2)} $\hat f_i$ is an $\epsilon_i$-Gromov-Hausdorff approximation (briefly GHA), $\epsilon_i\to 0$,

\noindent {\rm (1.13.3)} $\hat f_i$  is an $\epsilon_i$-Riemannian submersion i.e., for vector $\xi$ orthogonal to a $\hat F_i$, $e^{-\epsilon_i}\le \frac{|d\hat f_i(\xi_i)|}{|\xi_i|}\le e^{\epsilon_i}$.

\noindent {\rm (1.13.4)} The second fundamental form of $\hat F_i$, $|\op{II}(\hat F_i)|\le c(n)$.
\endproclaim

In this paper, a singular nilpotent fibration refers to a continuous map, $f_i: M_i\to X$, satisfying the commutative diagram (1.13.1), and $\tilde f_i$ may be a continuous map
satisfying (only) (1.13.2).

For applications of Theorem 1.13, see \cite{Ro2}.

\remark{Remark \rm 1.14} Let $f: M\to X$ be a singular nilpotent fibration,
and let $\pi: (\tilde M,\tilde p)\to (M,p)$ denote the Riemannian universal covering map.
The lifting singular fibration, $\tilde f: \tilde M\to \hat X$, is defined by that
a $\tilde f$-fiber is a component of $\pi^{-1}(f^{-1}(x)), x\in X$, and the following diagram commutes,
$$\CD \tilde F\to (F(\tilde M),\op{O(n)})@>\hat{\tilde f}>>(\hat Y,\op{O(n)})\\
@VV /\op{O(n)} V @VV /\op{O(n)}V\\
(\tilde M,\tilde p)@>\tilde f>>(\hat X=Y/\op{O(n)},\hat p)\endCD$$
From the homotopy exact sequence, $\tilde F\to F(\tilde M)\to \hat Y$, we conclude that $\pi_1(\tilde F)$ is abelian, thus $\tilde F$ is a nilpotent Lie group. Because $\hat Y$ is simply connected, the $\tilde F$-bundle is a principal $\tilde F$-bundle. Consequently, the lifting singular fibration on $\tilde M$ coincides with a $\tilde F$-action
on $\tilde M$ in the sense that a $\tilde F$-orbit coincides with a $\tilde f$-fiber (\cite{RY}).
\endremark

Let's now consider a sequence of collapsed $n$-manifolds of $\op{Ric}_{M_i}\ge -(n-1)$, and
$\forall\, x_i\in M_i$ is a $(\rho,v)$-LRVP (resp. a $(\delta,\rho)$-LRRP).

A first step in resolving (0.6.1) is to generalize Theorem 1.12 to the following.

\proclaim{Theorem 1.15} {\rm (Maximal collapsed manifolds with Ricci bounded covering geometry, \cite{HKRX}, \cite{Ro4})} Given $n, v>0$, there exists a constant $\epsilon(n,v)>0$, such that
if a compact $n$-manifold satisfies the following conditions, there is $\tilde p\in \tilde M$,
$$\op{Ric}_M\ge -(n-1),\quad \op{vol}(B_1(\tilde p))\ge v ,\quad \op{diam}(M)<\epsilon(n,v),$$
then $M$ is diffeomorphic to an infra-nilmanifold (in Theorem 1.12).
\endproclaim

Note that $\op{diam}(M)<\rho$ implies that $x\in M$ is a $(\rho,v)$-LRVP is equivalent to that
there is $p\in \tilde M$, $\op{vol}(B_1(\tilde p))\ge v$.

The proof of Theorem 1.15 in \cite{HKRX} relies on a smoothing method via Ricci flows (\cite{Ha}, \cite{Per2}): one first shows that points on $\tilde M$ are $(\delta,\rho)$-RP (Proposition 1.11), by which one applies Perel'man's pseudo-locality (\cite{Per2}) to
conclude that Ricci flows exists for a definite time which yields a metric satisfying Theorem 1.12, so Theorem 1.15 follows Theorem 1.12.

A second step is to prove Conjecture 0.6; Theorem 0.12 verifies (0.5.1)
in the case that $X$ is a manifold all whose points are $(\delta,\rho')$-RP (see Remark 1.17). Roughly speaking, the main part of Theorem 0.12 is a new proof of Theorem 1.15 without involving any smoothing technic (thus independent of Theorem 1.12) (\cite{Ro4}).

For convenience of readers, we will briefly describe the approach to Theorem 0.12 in \cite{Ro4}.
Recall that in \cite{Hu}, Huang constructed a fiber bundle map,
$F_i\to M_i@>f_i>>X$, by gluing locally defined $\delta$-splitting maps (\cite{CC1}). The non-degeneracy
for $\delta$-splitting maps and the glued map replies the transformation technique in
\cite{CJN} on a manifold of $(\rho,\delta)$-RP. Hence, the proof of Theorem 0.12 reduces to show
that $F_i$ is an infra-nilmanifold and the structural group can be reduced to an affine group that are unexpendable; which is unlikely to be verified due to a lack of regularity on $f_i$ (note that one may try a local smoothing technique (\cite{HW}); which likely fails due to that local smoothed metrics are not well-related).

To overcome the above problem, one has to find a direct proof of Theorem 1.15 without a smoothing. A direct proof is based on the following criterion for a compact manifold to be diffeomorphic to a nilpotent manifold.

\proclaim{Theorem 1.16} {\rm (A nilmanifold: a sequence of bundles, \cite{Ro4})}  A compact manifold $M_0$  is
diffeomorphic to a nilmanifold  if and only if $M_0$ carries a sequence of bundles:
$$M_1\to M_0\to B_1,\quad M_2\to M_0\to B_2,\quad \cdots, \quad M_{s+1}\to M_0
\to B_{s+1},$$
such that

\noindent {\rm (1.16.1)} Each $M_i$-fiber is a union of $M_j$-fibers for all $j\ge i$, $M_{s+1}=\{\op{pt}\}$.

\noindent {\rm (1.16.2)} For each $1\le i\le s$, the natural projection, $B_{i+1}\to B_i$, is a principal $T^{k_i}$-bundle ($k_1+\cdots +k_s=n$).
\endproclaim

Note that Theorem 1.16 is an alternative formulation of the criterion in \cite{Na} (\cite{Be}) that $M$ is diffeomorphic to a nilmanifold if and only if $M$ admits an iterated principal circle bundles:
$$S^1\to M\to M_1,\quad S^1\to M_1\to M_2,\quad \cdots \quad S^1\to M_n\to \op{pt}.$$

By a standard compactness argument (Theorem 1.4), we may consider a sequence in Theorem 1.15, $M_i@>\op{GH}>>\op{pt}$, and show that for $i$ large, $\tilde M_i$ is diffeomorhic to an infra-nilmanifold. Furthermore, with suitable rescaling and passing to a subsequence, we may assume that $\op{Ric}_{M_i}\ge -\epsilon_i(n-1)\to 0$, and show that $(\tilde M_i,\tilde p_i,\Gamma_i)@>\op{eqGH}>>(\tilde X,\tilde p,G)$ such that $\tilde X$ is isometric
to $\Bbb R^n$. Consequently, points in $\tilde M_i$ are $(\delta_i,\rho)$-RP,
$\delta_i\to 0$.

We start with Theorem 1.8, by which $M_i$ has a bounded normal covering space, $M_{i,0}$, with $\pi_1(M_{i,0})$ a nilpotent group. The proof consists of two steps: Step 1. (main work) $M_{i,0}$ satisfies the criterion of Theorem 1.16; Step 2. $M_i$ is diffeomorphic to an infra-nilmanifold. We point it out that in each step, the condition that $\tilde x_i\in \tilde M_i$ is a $(\delta_i,\rho)$-RP is essential.

To obtain a sequence of bundle structures on $M_{i,0}$, we blow-up $M_{i,0}$ at every
point $x_i\in M_{i,0}$ with an almost canonically determined rates, and show that the geometry and topology of a GH-limit is somewhat rigid i.e., independent of $x_i$ and a convergent subsequence. For instance, when a blow-up at $x_i$ by $\ell_{i,1}^{-1}$, $\ell_{i,1}$ is
the length of the first short generator for $\pi_1(M_i,x_i)$, then any GH-limit is a product flat torus $T^{k_1}$ with $\Bbb R^{n-k_1}$, $k$ is independent of $x_i$, and $\op{diam}(T^{k_1})\le d$, and $\op{injrad}(T^{k_1})\ge \epsilon>0$ ($\epsilon$ may depend on $\{M_i\}$). We then show that for $i$ large,
we obtain a locally finite open cover on $M_{i,0}$ by elementary T-structures in the sense of Cheeger-Gromov (\cite{CG1,2}). Employing the gluing method in \cite{CG2}, we construct a $T^{k_1}$-bundle on $M_{i,0}$, and show that, based on $\pi_1(M_{i,0})$ a torsion free nilpotent
group, this is a principal $T^{k_1}$-bundle.

Next, choosing the blow-up rate at $x_i$ by $\rho_{i,2}^{-1}$, $\ell_{i,2}$ is the length of
the first short generator that is not contained in $T^{k_1}$ at $x_i$, then the limit can be viewed as a blow up at a point in the orbit space of the principle $T^{k_1}$-bundle, $B_1$. Repeating the above operation, we conclude that $B_1$ admits a principal $T^{k_2}$-structure, thus $M_{i,0}\to B_1\to B_2$ is a fiber bundle and $B_1\to B_2$ is a principal $T^{k_2}$-bundle. Iterating the above operations, after at most $n$-steps, one finishes Step 1.

In the proof of Theorem 0.12, by performing the same successive blow ups at a point as in the above independent proof of Theorem 1.15, we obtain a locally finite open cover on $M$ with elementary N-structure, $f_\alpha: U_\alpha\to B_\alpha$, where $f_\alpha$ is a smooth trivial bundle map and $\Psi(\epsilon | n,v)$-GHA (a function in $\epsilon$ such that $\Psi(\epsilon|n)\to 0$ as $\epsilon\to 0$), with fiber an infra-nilmanifold, such that on overlaps, elementary N-structures satisfy the patching criterion in \cite{CFG}. Hence we are able glue elementary N-structures into a nilpotent fiber bundle with an affine structural group.

\remark{Remark \rm 1.17} In \cite{Ro4}, Theorem 0.12 was stated in the form that $X$ is a Riemannian manifold, and $x_i\in M_i$ is a $(\rho,v)$-LRVP (see (1.11.1)). In the proof, that $X$ is a Riemannian manifold is used only in two places: the construction of a local $(m,\delta)$-splitting map on $M_i$, and the gluing of local $(m,\delta)$-splitting map
via center of mass technique using a positive convexity radius of $X$.
It is easy to see that in both places, the condition that $X$ is a differentiable manifold
whose points are $(\delta,\rho)$-RP are enough in carrying out the above two
operations.
\endremark

\vskip4mm

\head 2. Proofs of Theorem 0.5, Theorem A, and Proposition 0.8
\endhead

\vskip4mm

\subhead 2.1. Proof of Theorem 0.5
\endsubhead

\vskip4mm

We first prove the following lemma, which is likely known to experts; and for a completeness we present a proof (cf. \cite{Wr}).

\proclaim{Lemma 2.1}  Let $\tilde M$ be an open contractible manifold. Let $G=T^k$ or $\Bbb Z_p^k$. If $G$ effectively acts on $\tilde M$, then the $G$-fixed point set is not empty and connected.
\endproclaim

In the proof, we will use the following Smith fixed point theorem.

\proclaim{Theorem 2.2} Let $M$ be a compact manifold, let $G=T^k$ or  $\Bbb Z_p^k$ ($p$ a prime) act effectively on $M$, and let $\ell$ be rationales $\Bbb Q$ for $G=T^k$, and $\ell=\Bbb Z_p$ for $G=\Bbb Z_p^k$. If $M$ is a homology $n$-sphere of coefficient $\ell$, then the $G$-fixed point set $F(M;G)$ is a homology sphere of coefficient $\ell$.
\endproclaim

\demo{Proof of Lemma 2.1}

Let $\tilde M^*$ denote the one-point compactification of $\tilde M$, and the $G$-action extends to $\tilde M^*$ and $G$ fixes $\tilde M^*\setminus \tilde M$. Because $\tilde M$ is contractible, $\tilde M^*$ is a homology $n$-sphere with a coefficient field $\ell$ in Theorem 2.2, thus the $G$-action has a non-empty fixed point set in $\tilde M$, otherwise the $G$-fixed point set in $\tilde M^*$ is a point, not a homology sphere; a contradiction to Theorem 2.2.

If the $G$-action has at least two fixed point components, $F_1, F_2$, then either the $G$-action on $\tilde M^*$ has two fixed point components, or the $G$-fixed point set in $\tilde M^*$ is not a homology sphere (e.g., the one point union of at least two homopology spheres), again a contradiction to Theorem 2.2.
\qed\enddemo

\proclaim{Lemma 2.3} Let the assumptions be as in Theorem 0.5. Then $\tilde M/G_0$ is not compact.
\endproclaim

\demo{Proof} Let $T^k$ denote the maximal compact subgroup of $G_0$. If $k=0$, then $G_0$ is
simply connected nilpotent Lie group acting freely on $\tilde M$, thus $\tilde M/G_0$ is contractible (not compact). Assume that $k\ge 1$. Because the $G_0$-action is not free (Lemma 2.1), we consider the free $G_0$-action on the principal $\op{Spin(n)}$-bundle over $\tilde M$, which is a double covering space of the oriented orthogonal frame bundle over $\tilde M$, with
the following commutative diagram,
$$\CD G_0\to (\op{Sp}(\tilde M),\op{Spin(n)})@> >> (\op{Sp}(\tilde M)/G_0,\op{Spin(n)})\\
@VV /\op{Spin(n)} V @VV /\op{Spin(n)}V\\
\tilde M@> >> \tilde M/G_0.
\endCD$$
Clearly, $\tilde M/G_0$ is not compact if and only if $\op{Sp}(\tilde M)/G_0$ is not compact.

Arguing by contradiction, assume that $\op{Sp}(\tilde M)/G_0$ is compact. Because $\tilde M$ is contractible, $\op{Sp}(\tilde M)$ is diffeomorphic to $\op{Spin(n)}\times \tilde M$, thus
$\op{Sp}(\tilde M)$ homotopically retracts to $\op{Spin(n)}$. We claim that $\op{Sp}(\tilde M)$ also homotopically retracts to a compact submanifold, $\op{Sp}(\tilde M)/(G_0\sim T^k)$,
which is a $T^k$-bundle over $\op{Sp}(\tilde M)/G_0$. Consequently, $\dim(\op{Spin(n)})=\dim(\op{Sp}(\tilde M)/(G_0\sim T^k))$ i.e., $n+k=\dim(G_0)$, a contradiction.

To verify the claim, observe that by Theorem 3.1, Exp: $T_{\bar e}(G_0/T^k)\to G_0/T^k$ is a diffeomorphism, thus the principal bundle, $G_0/T^k\to \op{Sp}(\tilde M)/T^k\to \op{Sp}(\tilde M)/G_0$, is equivalent to a vector bundle, $T_{\bar e}(G_0/T^k)\to\op{Sp}(\tilde M)/T^k\to \op{Sp}(\tilde M)/G_0$. Consequently, we may assume a homotopy retraction: $H_s: [0,1]\times \op{Sp}(\tilde M)/T^k\to \op{Sp}(\tilde M)/T^k$ such that $H_0=\op{id}_{\op{Sp}(\tilde M)/T^k}$ and $H_1: \op{Sp}(\tilde M)/T^k\to \op{Sp}(\tilde M)/G_0\subset \op{Sp}(\tilde M)/T^k$.
Consider the one-parameter family of pullback $T^k$-bundles by $H_s$:
$$\CD H_s^*(\op{Sp}(\tilde M))@>H_s^*>>\op{Sp}(\tilde M)\\
@VV /T^k V @VV /T^k V\\
[0,1]\times \op{Sp}(\tilde M)/T^k@>H_s>>\op{Sp}(\tilde M)/T^k,\endCD$$
Observe that each of the pullback $T^k$-bundle homotopy retracts to $T^k\to H_s^*(\op{Sp}(\tilde M))'\to \op{Sp}(\tilde M)/T^k$ as $[0,1]$ homotopy retracts to $0$,
and $H_s^*(\op{Sp}(\tilde M))'$ is one parameter family of $T^k$-bundles (isomorphic to
$T^k\to \op{Sp}(\tilde M)\to \op{Sp}(\tilde M)/T^k$). Because $H_s$ homotopy retracts
$\op{Sp}(\tilde M)/T^k$ to $\op{Sp}(\tilde M)/G_0$ (which is the `$0$' cross section), $H_1^*(\op{Sp}(\tilde M))'$ homotopy retracts to the pullback $T^k$-bundle, $\op{Sp}(\tilde M)/G_0\to \op{Sp}(\tilde M)/T^k$, denoted by $\op{Sp}(\tilde M)/(G_0\sim T^k)$ as in
the above; by now the claim has been verified.
\qed\enddemo

\demo{Proof of Theorem 0.5}

(0.5.1) Because the identity component $G_0$ is a connected nilpotent Lie group, its maximal compact subgroup is contained in the center i.e., a torus
$T^k$. We claim that $G_0$ has no compact subgroup i.e., $k=0$. Consequently,
$G_0$ acts freely on $\tilde M$.

Arguing by contradiction, assuming that $G_0$ has a maximal compact subgroup, a torus $T^k$ ($k\ge 1$). Because $\tilde M$ is contractible, by Lemma 2.1 the $T^k$-fixed point set, $F_0=F(T^k,\tilde M)\ne \emptyset$ is connected. Because $G_0<G$ is normal and $T^k<G_0$ is the unique maximal compact subgroup, $T^k$ is normal in $G$. Consequently, $G$ preserves $F_0$ i.e., $G(F_0)=F_0$, thus the discrete group, $G/G_0$, preserves the connected and closed subset of co-dimension $\ge 2$, $F_0/G_0\subset \tilde M/G_0$. Because $\tilde M$ is contractible, $\tilde M/G_0$ is not compact (Lemma 2.3). Let $\hat c(t)$ denote a ray at $\bar x_0\in F_0/G_0\subset \tilde M/G_0$ which is also normal to $F_0/G_0$. Because $\tilde M\to \tilde M/G_0$ is a submetry, there is horizontal ray $c(t)$ at $\tilde x_0\in F_0$ such that $\hat \pi(c(t))=\hat c(t)$. Observe that for any $\hat \gamma\in G/G_0$, $\hat \gamma(\hat c(t))$ is also a ray normal to $F_0/G_0$, which implies that $X=\hat X/(G/G_0)$ is not compact, a contradiction.

(0.5.2) Because $G_0$ has no compact subgroup, $G_0$ is diffeomorphic to $\Bbb R^s$. Because
$\tilde M$ is contractible, $\hat X=\tilde M/G_0$ is contractible. We claim that
$G/G_0$ acts freely on $\hat X$, thus $\hat X$ is a universal covering space of $X=\hat X/(G/G_0)$, an aspherical manifold.

Arguing by contradiction, assume a finite isotropy subgroup $\hat H<G/G_0$ that fixed $\hat x$. Without loss of generality, we may assume that $\hat H\cong \Bbb Z_p$, $p$ is a prime.
Then $\hat F_0=F(\hat X,\Bbb Z_p)\ne \emptyset$ is connected (Lemma 2.1). Because
any finite group of $G/G_0$ is virtually normal, the normalizer of $\hat H$, $N_{G/G_0}(\hat H)$, is a subgroup of $G/G_0$ with a finite index. Thus $N_{G/G_0}(\hat H)(\hat F_0)=\hat F_0$. By a similar argument as in the above, we conclude that $\hat X/N_{G/G_0}(\hat H)$ is not compact, a contradiction to that $X=\tilde X/G$ is compact.
\qed\enddemo

\subhead 2.2. Proof of Theorem A
\endsubhead

\vskip4mm

\demo{Proof of Theorem A}

By Theorem 1.13, we may assume that for $i$ large, there is a singular nilpotent fibration,
$f_i: M_i\to X$, and a $C^1$-close left invariant metric. Then the lifting singular
nilpotent fibration on $\tilde M_i$, $\tilde f_i: \tilde M_i\to \hat X$, coincides with the isometric-action of a connected nilpotent Lie group $\tilde F_i$-action i.e., a $\tilde f_i$-fiber is an $\tilde F_i$-orbit, thus the $\Gamma_i$-action preserves the
$\tilde F_i$-orbits fibers (see Remark 1.14).

By Theorem 0.5, $\tilde F_i$ contains no compact subgroup, thus $\tilde F_i$ is diffeomorphic to
$\Bbb R^k$. Let $\Lambda_i$ denote the subgroup that preserves a $\tilde F_i$-orbit. Because
the $\Gamma_i$-action preserves fibers of the principal $\tilde F_i$-bundle, $\Lambda_i<\Gamma_i$ is normal. Because $\Gamma_i$ is type $\Cal A_s$, any finite group of $\Gamma_i/\Lambda_i$ is virtually normal. Again by Theorem 0.5, the singular nilpotent fibration on $M_i$ is a fiber bundle, $F_i\to M_i@>f_i>>X$, thus $X$ is an aspherical differentiable manifold.

By Theorem 1.13, $X=Y/\op{O(n)}$, where $Y$ is $C^{1,\alpha}$-manifold on which $\op{O(n)}$
acts freely and isometrically. Consequently, the quotient metric on $X=Y/\op{O(n)}$
is defined by a $C^\alpha$-metric tensor.
\qed\enddemo

\vskip4mm

\subhead 2.3. Proof of Proposition 0.8
\endsubhead

\vskip4mm

\demo{Proof of Proposition 0.8}

Let $M_i@>\op{GH}>>X$ be as in Conjecture 0.7. Assume that for $i$ large, there is
a singular nilpotent fibration, $f_i: M_i\to X$, such that the induced metric on a $f_i$-fiber is almost left invariant. According to \cite{CFG}, there is a nearby left invariant metric whose
sectional curvature depends on $\max|\op{sec}_{M_i}|$ (obtained by averaging metric
by a (consistent) infinitesimal actions of a nilpotent Lie group). By now, the rest of the proof follows the proof of Theorem A,
and we conclude that $X$ is an aspherical manifold, differentialbe or homeomorphic, respectively.

It remains to show that in (0.6.1), all $x\in X$ is $(\delta,\rho')$-PR, $0<\rho'\le \hat \rho$.
By Theorem 0.9, for all $\tilde x_i\in \tilde M_i$, $\op{vol}(B_1(\tilde x_i))\ge v(n,d)>0$.
By (1.11.2), all points $\tilde x_i\in \tilde M_i$ are $(\delta,\tilde \rho)$-RP, thus all
points in $\tilde X$ are $(\delta,\tilde \rho)$-RP (Theorem 1.5), $0<\tilde \rho<\rho$.

Because $G_0$ acts freely and isometrically on $\tilde X$, $\tilde X\to \tilde X/G_0$ is submetry, it is clear that points in $\tilde X/G_0$ are $(\delta,\hat \rho)$-RP, $0<\hat \rho\le \tilde \rho$. Because discrete group $G/G_0$ acting freely and discontinuously on $\tilde X/G_0$ with a compact quotient $X$, points in $\tilde X/G$ are $(\delta,\rho')$-RP, for some $0<\rho'\le \hat \rho$.
\qed\enddemo

Note that in the above proofs of Theorem A and Proposition 0.8, Theorem 0.9 is not needed.

\vskip4mm

\head 3. Proofs of Theorem 0.11, Theorem B, and Proposition 0.13
\endhead

\vskip4mm

\subhead 3.1. Proof of Theorem 0.11
\endsubhead

\vskip4mm

In a proof of Theorem 0.11, our strategy is to introduce a left-invariant metric $\hat g_i$ on $\tilde M_i$ (identified with a simply connected nilpotent Lie group with identity $\tilde p_i$) and show that there is a constant $c(n,d)>0$ such that for any $R>0$,
$$B_R(\tilde p_i)\subset \bigcup_{\gamma_i\in \Gamma_i(R)}B_{c(n,d)}(\gamma_i(D(\tilde p_i,\hat g_i)))\subset B_{R+c(n,d)}(\tilde p_i),$$
where $D_i$ denote a Dirichlet fundamental domain at $\tilde p_i\in \tilde M_i$, $\Gamma_i(R)=\{\gamma_i\in \Gamma_i,\, \gamma_i(D_i)\cap B_R(\tilde p_i)\ne \emptyset\}$.

Given a higher homotopy class of $\tilde X$ represented by a continuous map, $\phi: S^k\to \tilde X$, $\phi(S^k)\subset B_R(\tilde p)$ for some $R>0$, we shall show that there is homotopy
$H_s'$ deforming $f_i^{-1}(\phi(S^k_1))$ to $\tilde p_i$ inside
$B_{R+c(n,d)}(\tilde p_i)$ ($0\le s\le 1$), where for $i$ large $f_i:  B_{R+c(n,d)}(\tilde p_i)\to B_{R+c(n,d)}(\tilde p)$ is a diffeomorphism (Theorem 1.5) or homeomorphism (Theorem 1.6). Then $f_i\circ H_s'\circ f_i^{-1}$ deforms $\phi(S^k_1)$ to $\tilde p$, thus i.e., $\tilde X$ is contractible.

Our construction of $H_s'$ relies on the following result.

\proclaim{Theorem 3.1} {\rm (\cite{Sa1,2})} If $G$ is an exponential Lie group, then $G$ is simply connected and solvable, but the converse fails to hold. Moreover, a simply connected nilpotent Lie group is exponential.
\endproclaim

We now define a left-invariant metric on $\tilde M_i$ based on information from a set of (Gromov) short generators for $\Gamma_i=\pi_1(M_i,p_i)$ as follows.

Recall that for a compact $n$-manifold $M$ of $\op{diam}(M)=d$, $p\in M$,
$\Gamma=\pi_1(M,p)$ can be generated by a finite number of elements whose displacement at
$\tilde p$ is $\le 2d+\epsilon\le 3d$ (for any $0<\epsilon<\op{injrad}(M)$, the injectivity radius). A set of short generators for $\pi_1(M,p)$ is picked up in order as follows: let $\gamma_1\in \Gamma$ such that $|\gamma_1|=d(\tilde p,\gamma(\tilde p))$ achieves the minimum among all $\gamma\in \Gamma$. Let $\gamma_2\in \Gamma\setminus \left<\gamma_1\right>$ (the subgroup generated by $\gamma_1$) such that
$|\gamma_2|$ is the minimum among all $\gamma\in \Gamma\setminus \left<\gamma_1\right>$. Iterating this process, one gets a set of short generators, denoted by $S_{\tilde p}(\Gamma)$.

For a finitely generated torsion free nilpotent group $\Lambda$, we call the minimal number of
generators the rank of $\Lambda$.

\proclaim{Lemma 3.2} Let a sequence of compact $n$-manifolds $M_i@>\op{GH}>>X$ satisfying
$$\CD (\tilde M_i,\tilde p_i,\Gamma_i)@>\op{eqGH}>>(\tilde X,\tilde p,G)\\
@VV \pi_i V @VV /GV \\
(M_i,p_i)@>\op{GH}>>(X=\tilde X/G,p),\endCD,\quad \cases \op{Ric}_{M_i}\ge -(n-1)\\ 0<\op{diam}(X)=d<\infty.\endcases$$
If $\Gamma_i(\epsilon)$ is a torsion free nilpotent group of rank $s$, then $G_0$ is diffeomorphic to $\Bbb R^a\times T^b$, $a\ge s, b\ge 0$.
\endproclaim

\demo{Proof} We proceed by induction on $s$, starting with $s=1$. Because $\Gamma_i(\epsilon)$ is torsion free nilpotent group of rank one, $\Gamma_i(\epsilon)\cong \Bbb Z$. Because
$\Gamma_i(\epsilon)(\tilde p_i)$ is not bounded, which converges $G_0(\tilde p)$ which is not bounded i.e., $G_0(\tilde p)$ is an open submanifold of dimension $\ge 1$. Consequently,
$G_0$ is a non-compact nilpotent Lie group, thus $G_0\cong \Bbb R^a\times T^k$, $a\ge 1$ and $k\ge 0$, where $T^k$ is a maximal compact subgroup.

Consider torsion free nilpotent $\Gamma_i(\epsilon)$ of rank $s+1$. Then $\Gamma_i(\epsilon)$ has a non-trivial center, let $\gamma_i\in C(\Gamma_i(\epsilon))$ such that $|\gamma_i|\to 0$. Let $\hat M_i=\tilde M_i/\left<\gamma_i\right>, \hat \Gamma_i(\epsilon)=\Gamma_i(\epsilon)/\left<\gamma_i\right>$,
$\left<\gamma_i\right>\to G_1$, which is contained in the center of $G_0$. Consider
the following commutative diagram commutes (see (1.3.3)):
$$\CD (\tilde M_i,\tilde p_i,\Gamma_i(\epsilon))@>\op{eqGH}>>(\tilde X,\tilde p,G_0)\\
@VV /\left<\gamma_i\right> V @VV /G_1 V \\
(\hat M_i,\hat p_i, \hat \Gamma_i(\epsilon))@>\op{GH}>>(\hat X=\tilde X/G_1,\hat p, G_0/G_1),\endCD$$
Because $\hat \Gamma_i(\epsilon)$ is torsion free nilpotent group of rank $s$ and $\hat \Gamma_i(\epsilon)\to G_0/G_1$, by induction we have that $G_0/G_1\cong \Bbb R^a\times T^{k'}, a\ge s, b\ge 0$ (no that in the induction, that $\tilde M_i$ is simply connected is not required). Because $G_1=\Bbb R^{a'}\times T^{b'}$, $a'\ge 1, b'\ge 0$, we conclude that $G_0\cong \Bbb R^{a+a'}\times T^b, a+a'\ge s+1, b\ge 0$.
\qed\enddemo

Identifying $\tilde M_i$ with a simply connected nilpotent Lie group with identity $e=\tilde p_i$, one may define an inner product on $T_{\tilde p}\tilde M_i$, by assigning an inner product for any linearly independent $n$ vectors in $T_{\tilde p_i}\tilde M_i$. Note that because $\Gamma_i<\tilde M_i$ acts on $\tilde M_i$ isometrically with respect to the left-invariant metric, the left-invariant metric descends to a metric on $M_i$, denoted by $\hat g_i^*$.

First, each short generator, $\gamma_{i\alpha}\in S_{\tilde p_i}(\Gamma_i)$, determines a vector in $T_{\tilde p_i}\tilde M_i$:
$$u_{i\alpha}=d_i(\tilde p_i,\gamma_{i\alpha}(\tilde p_i))\uparrow_{\tilde p_i}^{\gamma_{i\alpha}(\tilde p_i)},$$
where $\uparrow_{\tilde p_i}^{\gamma_{i\alpha}(\tilde p_i)}$ denotes the direction of a minimal geodesic from $\tilde p_i$ to $\gamma_{i\alpha}(\tilde p_i)$ (which may not be unique). Then the set of vectors, $\{u_{i\alpha}, \,\, \gamma_{i\alpha}\in S_{\tilde p_i}(\Gamma_i)\}$, spans $T_{\tilde p_i}\tilde M_i$. We will start with $u_{ik}=u_{in}$ (the last one) and extend it linearly independently, in the reverse order, to a basis for $T_{\tilde p_i}\tilde M_i$, denoted by $\{u_{i1},...,u_{in}\}$. Note that because $\left<\gamma_{i1},...,\gamma_{in}\right><\Gamma_i$ is a torsion free nilpotent group, $\gamma_{i1}$ is from the center of $\Gamma_i$.

We now assign an inner product for the basis, $\{u_{i1},...,u_{in}\}$,
$$\hat g_i(u_{i\alpha},u_{i\beta})=d_i(\tilde p_i,\gamma_{i\alpha}(\tilde p_i))\cdot d_i(\tilde p_i,\gamma_{i\beta}(\tilde p_i))\cos\measuredangle (\uparrow _{\tilde p_i}^{\gamma_{i\alpha}(\tilde p_i)}\uparrow _{\tilde p_i}^{\gamma_{i\beta}(\tilde p_i)}).$$

Let $d_i=\op{diam}(M_i)$ and $\hat d_i=\op{diam}(M_i,\hat g_i^*)$.

\proclaim{Lemma 3.3} {\rm (Ratio of diameters estimate)} Let the assumptions be as in Proposition 0.8, and let $\hat  g_i^*$ be the left-invariant metric on $M_i$ defined in the above. Then there exist constants, $a(n), b(n,d)>0$, such that for $i$ large,
$$a(n)\hat d_i\le d_i\le b(n,d)\hat d_i.$$
\endproclaim

\demo{Proof} We first claim that there is a constant $c(n)$ such that the following inequality holds:
$$\ell_i<\hat d_i\le c(n)\ell_i,\quad \ell_i=\max\{|\gamma_{i,j}|,\,\gamma_{i,j}\in S_{\tilde p_i}(\Gamma_i)\}.\tag 3.3.1$$
We proceed a proof by induction on $n$, starting with $n=2$: $\hat g_i$ is flat,
and diameter $\hat d_i\le \sqrt 2\ell_i$. Assume the claim holds for $\dim(M_i)=n-1$.

In general, by the choice of $\gamma_{i1}$ we may assume that $\left<\gamma_{i1}\right>$ is contained in the center of $\Gamma_i$. Then there is a $\Bbb R^1$-subgroup of $\tilde M_i$
that contains $\left<\gamma_{i1}\right>$. Let $\hat M_i=(\tilde M_i/\Bbb R^1/(\Gamma_i/\left<\gamma_{i1}\right>$, which is a nilmanifold of dimension $n-1$. Applying the inductive assumption on $((\tilde M_i,\hat g_i)/\Bbb R^1)/(\Gamma_i/\left<\gamma_{i1}\right>)$, $\hat d_i\le c(n-1)\ell_i^*\le c(n-1)\ell_i$ that $S^1\to M_i\to \hat M_i$ is a Riemannian submersion ($S^1=\Bbb R^1/\left<\gamma_{i1}\right>$), it is clear that $\hat d_i\le c(n-1)\ell_i+\ell_i=(c(n-1)+1)\ell_i=c(n)\ell_i$.

We then claim that there is a constant $b(n,d)>0$ such that the following inequality holds:
$$\frac 13\ell_i\le d_i\le b(n,d)\ell_i.\tag 3.3.2$$
Note that the desired inequality follows from a combination of (3.3.1) and (3.3.2), with
$a(n)=\frac 1{3c(n)}$.

Because $\Gamma_i$ is generated by elements in $S_{\tilde p_i}(\Gamma_i)$, whose distortions at $\tilde p_i$ is less than $3d_i$, $\frac 13\ell_i\le d_i\to d>0$. We now prove the other inequality by a contradiction: assuming, passing to a subsequence, that $\ell_i\to 0$ as $i\to \infty$. By Lemma 3.2, $G_0\cong \Bbb R^n\times T^k, k\ge 0$, thus $G_0$ acts transitively on $\tilde X$; so $X=\tilde X/G_0$ is a point, a contradiction to $\op{diam}(X)=d>0$.
\qed\enddemo

\demo{Proof of Theorem 0.11}

Let $M_i@>\op{GH}>>X$ be as in Theorem B, and let $\tilde X$ be as in (0.10). Given a higher homotopy class of $\tilde X$ represented by a continuous map, $\phi: S^k\to \tilde X$, $\phi(S^k)\subset B_R(\tilde p)$ for some $R>0$.

Let $D(\tilde p_i)$ denote a Dirichlet fundamental domain at $\tilde p_i$ (with respect to $\tilde g_i$). Let $\Gamma_i(R)=\{\gamma_i\in \Gamma_i, \, \gamma_i(D_i(\tilde p_i))\cap B_R(\tilde p_i)\ne \emptyset\}$ (a finite set). Then
$$B_R(\tilde p_i)\subset \bigcup_{\gamma_i\in \Gamma_i(R)}\gamma_i(D(\tilde p_i)).$$
By Lemma 3.3, we get that $\gamma_i(D(\tilde p_i))\subset B_{c(n,d)}(\gamma_i(\tilde p_i))$, $c_1(n,d)=b(n,d)+3c(n)$. Because $d(\tilde p_i,\gamma_i(\tilde p_i))\le R+d$, we derive
$$\bigcup_{\gamma_i\in \Gamma_i(R)}\gamma_i(D(\tilde p_i))\subset \bigcup_{\gamma_i\in \Gamma_i(R)}B_{c(n,d)}(\gamma_i(D(\tilde p_i,\hat g_i)))\subset B_{R+c(n,d)}(\tilde p_i),$$
where $c(n,d)=d+c_1(n,d)$. For $i$ large, we may assume a diffeomorphism (Theorem 1.5) or a homeomorphism (Theorem 1.6),
$$f_i: B_{R+c(n,d)}(\tilde p_i)\to B_{R+c(n,d)}(\tilde p).$$
Let $H_s$ denote the radial contraction from $B_{R+c(n,d)}(0_i)\subset T_{p_i}\tilde M_i$ to $0_i$. Identify $(\tilde M_i,\tilde p)$ with a simply connected nilpotent Lie group (diffeomorphically or homeomorphically),
because $\op{Exp}: T_{\tilde p_i}\tilde M_i\to \tilde M_i$ is a diffeomorphism or a homeomorphic (if
$M_i$ is homeomorphic to a nilmanifold),
$$\exp_{\tilde p_i}(H_s(\op{Exp}^{-1}(f_i^{-1}(\phi(S^k_1)))))\subset B_{R+c(n,d)}(p_i),\quad 0<s\le 1,$$ defines a deformation contraction, $f_i^{-1}(\phi(S^k_1))$ to $\tilde p_i$, inside $B_{R+c(n,d)}(\tilde p_i)$, thus $f_i\circ \exp_{\tilde p_i}\circ H_s\circ \op{Exp}^{-1}\circ f_i^{-1}: [0,1]\times \phi(S^k_1)\to \tilde X$ defines a deformation contraction of $\phi(S^k_1)$ to $\tilde p$, i.e.,
$\tilde X$ is contractible.
\qed\enddemo

\remark {Remark \rm 3.4} Inspecting the above proof, it is clear that Theorem 0.11 is valid when replacing $\tilde M$ (identified as a simply connected nilpotent Lie group) by
an exponential solvable Lie group.
\endremark

\vskip4mm

\subhead 3.2. Proofs of Theorem B and Proposition 0.13
\endsubhead

\vskip4mm

Let $M_i@>\op{GH}>>X$ be as in Theorem B, and let $(\tilde X,G)$ be in (0.10). We first extend Theorem 0.5 to ($\tilde X,G)$, where $\tilde X$ is a contractible differentiable, or a topological $n$-manifold, respectively (Theorem 0.11).

\proclaim{Lemma 3.5} Let $M_i@>\op{GH}>>X$ be as in Theorem B, and let $(\tilde X,G)$ be as in (0.10). Then $G$ acts freely on $\tilde X$, thus $\tilde X/G$ is a differentiable, or topological
manifold, respectively.
\endproclaim

Note that in (B2), $\tilde X$ is a smoothable Alexandrov space of curvature $\ge -1$
i.e., iterated space of directions at a point are spheres (\cite{Ka1}). According to
Proposition 2.7 in \cite{BNS}, if a closed group of isometries $G$ freely acts on
a smoothable Alexandrov space $X$, then $X/G$ is a topological manifold.

\demo{Proof of Lemma 3.5}

Let $\op{Sp}(\tilde M_i)$ denote the $\op{Spin(n)}$ bundle over $\tilde M_i$ equipped with a canonical metric. Passing to a subsequence, we may assume the following commutative diagram,
$$\CD (\op{Sp}(\tilde M_i),(\tilde p_i,u_i),\op{Spin(n)},\Gamma_i(\epsilon),\Gamma_i)@>\op{eqGH}>>(\op{Sp}(\tilde X),(\tilde p,u),\op{Spin(n)},G_0<G)\\
@VV /\op{Spin(n)}V @VV /\op{Spin(n)}V\\
(\tilde M_i,\tilde p_i,\Gamma_i(\epsilon),\Gamma_i)@>\op{eqGH}>>(\tilde X,\tilde p,G_0<G)
\endCD$$
By Theorem 0.11, $\tilde X$ is contractible $n$-manifold, and because $\Gamma_i$ is nilpotent,
$G$ is nilpotent. Note that the Ricci curvature of the canonical metric on $\op{Sp}(\tilde M_i)$ may not have a uniform lower bound; so the equivariannt GH-convergence is from that for any $R>0$, there is a uniform upper bound on the number of an $\epsilon$-net in $B_R(\tilde p_i,u_i)$, $s(\epsilon)\cdot \ell_i(R,\epsilon)$, where $s(\epsilon)$ denotes an $\epsilon$-net in $\op{Spin(n)}$ with a bi-invariant metric, and $\ell_i(R,\epsilon)$ is the
number of $\epsilon$-net in $B_R(\tilde p_i,\tilde M_i)$.

Note that one may view the free isometric $G$-action on $\op{Sp}(\tilde X)$ (equipped with a
length metric) as the lifting of isometric $G$-action on $\tilde X$. Consequently, $\tilde X/G_0$ is not compact. By now the rest of proof follows from the proof of Theorem 0.5.
\qed\enddemo

\demo{Proof of (B1)}

By Theorems 0.9, 1.5 and 0.11, in (0.10) $\tilde X$ is a contractible $n$-manifold. Because
$M_i$ is nilmanifold, $G$ is nilpotent. 
Applying Lemma 3.5 to $(\tilde X,G)$, we conclude that $X$ is an aspherical manifold, and by Proposition 0.7, points in $X$ are $(\delta,\rho')$-RP.

We now apply Theorem 0.12 to $M_i@>\op{GH}>>X$, and conclude that for $i$ large there is a smooth fiber bundle map, $F_i\to M_i@>f_i>>X$, with $F_i$ an infra-nilmanifold and affine structural group.

Let $\pi_i: (\tilde M_i,\tilde p_i)\to (M_i,p_i)$ denote the Riemannian universal covering map, and consider the lifting bundle, $\tilde F_i\to \tilde M_i\to \hat X$, where $\tilde F_i$ is a component of $\pi^{-1}_i(F_i)$. Because $\tilde M_i$ is simply connected, $\hat X$ is simply connected, and $\hat \pi: (\hat X,\hat p)\to (X,p)$ is
a universal covering map. Because $\hat X$ is contractible, $\tilde F_i\to \tilde M_i\to \hat X$
is a trivial bundle, $\tilde M_i\cong \tilde F_i\times \hat X$, thus $\tilde F_i$ is simply connected i.e., $\tilde F_i$ is a simply connected nilpotent Lie group.

We now identify $(\tilde M_i,\tilde p_i)$ with a simply connected nilpotent Lie group of $e=\tilde p_i$, and $\Gamma_i$ a co-compact lattice, and identify $(\tilde F_i,\tilde p_i)$ with a simply connected nilpotent Lie group of $e=\tilde p_i$, and a co-compact lattice $\Gamma_i(\epsilon)=\pi_1(F_i)$ (see (1.10.2)).

Observe that the inclusion map, $F_i\hookrightarrow M_i$, induces a homomorphic embedding,
$\Gamma_i(\epsilon)\hookrightarrow \Gamma_i$, which, by Malc\'ev rigidity, induces an isomorphic embedding on simply connected nilpotent Lie groups, $\tilde F_i\hookrightarrow \tilde M_i$. By now the principal $\tilde F_i$-bundle on $\tilde M_i$ is identified as a normal subgroup Lie group $\tilde F_i$-action of the nilpotent Lie group $\tilde M_i$, thus
$\hat X$ is diffeomorphic to the simply connected nilpotent Lie group, $\tilde M_i/\tilde F_i$, with a co-compact lattice $\Gamma_i/\Gamma_i(\epsilon)$ (see (1.10.2)). By now we conclude that $X$ is diffeomorphic to $\hat X/(\Gamma_i/\Gamma_i(\epsilon))$, a nilmanifold.
\qed\enddemo

In the proof of (B2), we shall apply the following Farrell-Hsiang's result which confirms the Borel conjecture for nilmanifolds (the Borel conjecture asserts: two homotopy equivalent compact aspherical manifolds are homeomorphic).

\proclaim{Theorem 3.6} {\rm (\cite{FH})} Two homotopy equivalent compact nilmanifolds are homeomorphic.
\endproclaim

\demo{Proof of (B2)}

By Theorems 0.9, and 1.6 and 0.11, $\tilde X$ is a contractible topological $n$-manifold. Because
$M_i$ is nilmanifold, $G$ is nilpotent. 
By Lemma 3.5, $X$ is a topological aspherical manifold of $\pi_1(X)\cong G/G_0$, a nilpotent
group. By the homotopy theory, $X$ is homotopically equivalent to a nilmanifold. Then (B2) follows from Theorem 3.6 i.e., $X$ is homeomorphic to a nilmanifold.
\qed\enddemo

\remark{Remark \rm 3.6} The above proof of (B2), restricting to the case that $M_i$ is homeomorphic to a torus, is different from the proof of Theorem 0.2 in \cite{BNS}.
\endremark

\demo{Proof of Proposition 0.13}

Consider a sequence of aspherical $n$-manifolds, $M_i@>\op{GH}>>X$, satisfying
the conditions in Conjecture 0.7. By Theorem 0.8, 1.5 and 1.6, $\tilde X$ in (0.10) is an $n$-dimensional manifold, differentiable or homeomorphic, respectively. Assume that $\tilde X$  is contractible (Conjecture 0.13). By (1.10.1) $G_0$ is nilpotent, and by (1.10.2) for $i$ large $G/G_0$ is isomorphic to $\Gamma_i/\Gamma_i(\epsilon)$. Because $\Gamma_i$ is type $\Cal A_s$, any finite subgroup of $G/G_0$ is virtually normal. We now apply Lemma 3.5 to conclude that $X=\tilde X/G$ is an aspherical manifold, differentiable or topological, respectively i.e., Conjecture 0.7 is proved.
\qed\enddemo

\vskip4mm

\head 4. Proof of Theorem C
\endhead

\vskip4mm

\demo{Proof of Theorem C}

Let $M_i@>\op{GH}>>X$ be a sequence of compact aspherical $3$-manifolds satisfying
$$\op{Ric}_{M_i}\ge -1,\quad 0<\op{diam}(X)=d<\infty,\quad \op{vol}(M_i)\to 0,$$
and the commutative diagram (0.9), $X=\tilde X/G$.

By Theorem 0.9, $\op{Vol}(B_1(\tilde p_i))\ge v>0$, and by Theorem 1.7 in \cite{Si}, $\tilde X$ is an open $3$-manifold. Because $0<\op{diam}(X)=d<\infty$, $1\le \dim(G)\le 2$.

Case 1. Assume $\dim(G)=1$. By Lemma 3.2, $G_0=\Bbb R^1$. Because $\Bbb R^1$ acts freely on $\tilde X$, $\tilde X/G_0$ is open surface. We claim that $\tilde X/G_0$ is simply connected, thus $\tilde X/G_0$ is homeomorphic to $\Bbb R^2$. By Lemma 3.5, $G/G_0$ acts freely
on $\tilde X/G_0$, thus $X$ is an aspherical manifold.

If $\tilde X/G_0$ is not simply connected, by (1.10.4) $(\tilde X/G_0)/H$ is simply connected open surface, where $H$ is generated by finite isotropy groups of the $G/G_0$-action on $\tilde X/G_0$. Because the fixed point set of each isotropy group is isolated, that $\tilde X/G_0$ is not simply connected implies that $(\tilde X/G_0)/H$ not simply connected, a contradiction.

Case 2. Assume that $\dim(G)=2$. Then $G_0=\Bbb R^1\times S^1$ or $\Bbb R^2$.
Because $\tilde X/G_0$ is not compact, $\tilde X/G_0=\Bbb R$ or $\Bbb R_+$. Note that $\Bbb R_+$ has a trivial isometry group, thus $G/G_0=e$, a contradiction to that $X$ is compact. If $\tilde X/G_0=\Bbb R$, then that $\pi_1(M_i)$ is type $\Cal A_s$ implies that $G/G_0$ has no torsion element (otherwise, $G/G_0=\Bbb Z_2\rtimes \Bbb Z$, a contradiction). Consequently, $X=\tilde X/G$ is a circle, an aspherical $1$-manifold.
\qed\enddemo

Next, we supply a justification of the following assertion in discussion before Definition 0.4.

\proclaim{Lemma 4.1} There is no compact simply connected $n$-manifold $M$ ($n\le 4$)
which admits and only admits a free $T^1$-action.
\endproclaim

\demo{Proof} For $n=2$, $M$ is a sphere and any $T^1$-action on $M$ has a fixed point.

For $n=3$, if $M$ admits a free $T^1$-action, then $M/T^1$ is a compact simply
connected surface, thus $M/T^1$ is homeomorphic to $S^2$, and therefore
$M$ is homeomorphic to $S^3$ (one may conclude this directly from the Pincar\'e conjecture).
But $S^3$ also admits not free $T^1$ and $T^2$-actions.

For $n=4$, then the Euler characteristic of $M$ is $\ge 2$, thus $M$ admits no
free $T^1$-action.
\qed\enddemo

The following example was mentioned in Introduction following Conjecture 0.6.

\example{Example 4.2} Let $M=S^k\times S^{n-k}$, $k\ge 2$ and $n\ge 5$. In \cite{Ot},
Otsu constructed on a family of metrics $g_i$ on $M$, such that $\op{Ric}_{g_i}\ge n-1$, $\op{vol}(M,g_i)\ge v>0$ and $\op{diam}(g_i)\to \pi$. The following
property on $S^3\times S^3$ was pointed out to the author by J. Pan: $g_i$ has a free
free $S^1$-symmetry, $(S^3\times S^3,g_i)@>\op{GH}>>\tilde X$, where $\tilde X$ is
a spherical suspension over $S^3\times S^2$ equipped with a metric of diameter $\pi$.
Let $\Bbb Z_i<S^1$. Then
$$\CD (S^3\times S^3,g_i,\Bbb Z_i)@>\op{eqGH}>>(\tilde X,G=S^1)\\
@VV /\Bbb Z_i V @VV /S^1V\\
(S^3\times S^3)/\Bbb Z_i@>\op{GH}>>X=\tilde X/S^1.
\endCD$$
Because the $S^1$-action fixes a vertex of $\tilde X$, for all $i$ large, $(S^3\times S^3)/\Bbb Z_i$ admits no singular nilpotent fibration, because otherwise the lifting of the singular nilpotent fibration on $S^3\times S^3$ is a $S^1$-action without fixed points, which converges to the $S^1$-action on $\tilde X$, which must have no fixed point, a contradiction.
\endexample

In Example 4.2, $(S^3\times S^3)/\Bbb Z_i$ admits an isometric $S^1$-action with
isolated fixed points; which may be treated as a generalized singular fibration
if one drops the condition that each fiber is an infra-nilmanifold of positive dimension.

To partially justify the existence of a generalized singular nilpotent fiber bundle on a collapsed manifold $M$ of $\op{Ric}_M\ge -(n-1)$ and $x\in M$ is a $(\rho,v)$-LRVP, let first
recall the following property on a Ricci limit space of a collapsing sequence of $n$-manifolds
of $(\rho,v)$-LRVP (\cite{Ro3}).

\proclaim{Lemma 4.3} {\rm (A Ricci limit space of a collapsing sequence with LRVP)} Let a sequence of complete $n$-manifolds, $(M_i,p_i)@>\op{GH}>>(X,p)$, satisfying
$$\op{Ric}_{M_i}\ge -(n-1),\quad \forall\, \text{ $x_i\in M_i$ is a $(\rho,v)$-LRVP}.$$

\noindent {\rm (4.3.1)} For $x\in X$, any tangent cone at $x$ is a metric cone of dimension $m\in \Bbb N$, thus the Hausdorff dimension $\dim_H(X)=m$.

\noindent {\rm (4.3.2)} The set of $\epsilon$-regular points in $X$, $X_\epsilon$, is an open $m$-manifold which is dense in $X$.

\noindent {\rm (4.3.3)} For small $0<\rho<1$, points in $X_0=\{x\in X_\epsilon, \, d(x, X\setminus X_\epsilon)\ge \rho\}$ are $(\delta,\frac \rho4)$-RP, $0<\delta\le \delta(n)$.
\endproclaim

Because \cite{Ro3} is in chinese, for convenience of readers we present the following proof.

\demo{Proof} (4.3.1) Let $x\in X$ and let $r_j\to \infty$ such that $(r_jX,x)@>\op{GH}>>(C_xX,x)$. We may assume $x_i\in M_i$, such that
$x_i\to x$. We start with the following commutative equivariant convergence (\cite{Xu}),
$$\CD (\widetilde{B_\rho(x_i)},\tilde x_i,\Gamma_i)@>\op{eqGH}>>(\tilde X, \tilde x,G) \\@VV \pi_i V @VV /G V\\
(B_\rho(x_i),x_i)@>\op{GH}>>(B_\rho(x)=\tilde X/G,x),
\endCD$$
where $\pi_i: (\widetilde {B_\rho(x_{i_j})},\tilde x_{i_j})\to (B_\rho(x_{i_j}),x_{i_j})$ denotes the Riemannian universal covering map.
Then passing to a subsequence, $(r_jB_\rho(x_{i_j}),x_{i_j})@>\op{GH}>>(C_xX,x)$,
and the following commutative diagram:
$$\CD (r_jB_\rho(\tilde x_{i_j}),\tilde x_{i_j},\Gamma_{i_j})@>\op{eqGH}>>(C_{\tilde x}\tilde X, \tilde x,dG) \\@VV \pi_i V @VV /G V\\
(r_jB_\rho(x_{i_j}),x_{i_j})@>\op{GH}>>(C_xX=C_{\tilde x}\tilde X/dG,x),
\endCD$$
where $dG$ denotes the `differential' of the $G$-action at $\tilde x$ i.e., the limit of
the subgroup of $G$ generated by $G_0$ (the identity component) and isotropy group at $\tilde x$. Because $\op{vol}(B_\rho(x_{i_j}))\ge v>0$, $C_{\tilde x}\tilde X=\Bbb R^k\times C(\Sigma_{\tilde x})$ is a metric cone of dimension $n$, and $\op{diam}(\Sigma_{\tilde x})<\pi$ (\cite{CC1}). We may assume that $dG=\Bbb R^l\times G_1$ ($l\le k$), $\Bbb R^l$ acts by translattions on the $\Bbb R^k$-factor and $G_1$ is a compact subgroup acting on the cross section of $C_{\tilde x}X/\Bbb R^l=C(\Sigma_{\tilde x}(\tilde X/\Bbb R^l))$. By now it is clear that $C_xX \cong C(\Sigma_{\bar x}(C_{\tilde x}(\tilde X/\Bbb R^l)/G_1))$, a metric cone.

(4.3.2) This is a consequence of a result in \cite{CC3}: if a Ricci limit space satisfies that
every tangent cone is a pole space, then $\dim_H(X)\in \Bbb N$.

(4.3.3) Let $x\in X_0$, and we may assume $x_i\in M_i$ such that $x_i\to x$. Because $B_{\frac \rho2}(x)\subset X_\epsilon$, consider the following commutative equivariant GH-convergence,
$$\CD (\widetilde{B_{\frac \rho2}(x_i)},\tilde x_i,\Gamma_i)@>\op{eqGH}>>(\tilde X, \tilde x,G) \\@VV \pi_i V @VV \op{proj}V\\
(B_{\frac \rho2}(x_i),x_i)@>\op{GH}>>(B_{\frac \rho2}(x)=\tilde X/G,x).
\endCD$$
Because $\tilde X$ is an $n$-manifold and $B_{\frac \rho2}(x)$ is a manifold, the isometric
(Lie group) $G$-action on $\op{proj}^{-1}(B_{\frac \rho2}(x))$ is free. Because $\tilde x$ is $(\delta,\frac \rho2)$-RP ((1.11.2)), $x$ is a $(\delta,\rho')$-RP, $0<\rho'<\frac \rho2$.
\qed\enddemo

Combining Lemma 4.3 and Theorem 0.12, in the case that $0<\op{diam}(X)=d<\infty$ we conclude
a fiber bundle, $F_i\to f_i^{-1}(X_0)\to X_0$, with $F_i$ an infra-nilmanifold and affine
structural group.

\example{Open problem 4.4} (Generalized singular nilpotent fibration, \cite{HRW}) Let a sequence of compact $n$-manifolds, $M_i@>\op{GH}>>X$, such that
$$\op{Ric}_{M_i}\ge -(n-1),\quad 0<\op{diam}(X)=d<\infty,\quad x_i\in M_i \text{ is a $(\rho,v)$-LRVP}.$$
Then for $i$ large, there is continuous $\epsilon_i$-GHA, $f_i: M_i\to X$,
such that $F_i\to f_i^{-1}(X_0)\to X_0$, is smooth fiber bundle.
\endexample

Note that a positive answer to Open Problem 4.4 would likely provide information on a singular $f_i$-fiber; which could be a point (Example 4.2), or whether which is a topological infra-nilamifold?

\vskip4mm

\head Appendix
\endhead

\vskip4mm

In the appendix, we will supply a proof (1.3.3) due to J. Pan.

\demo{Proof of (1.3.3)}

Let $(X,p,\Lambda,<\Gamma):=(X_i,p_i,\Lambda_i<\Gamma_i)$, let $\pi: \Gamma\to \frac \Gamma\Lambda$, and let $\bar X=X/\Lambda$ denote the quotient space with quotient metric $d_{\bar X}$. Following a standard argument, one obtains (1.3.3) if for $R>0$, the following
equality holds; \footnote{In Reamrk 3.5 in \cite{FY}, it was claimed that in general the
equality does not hold.}
$$\frac \Gamma\Lambda(R)=\frac{\Gamma(R)\Lambda}\Lambda.$$
where $\Gamma(R)=\{\gamma\in \Gamma, d_X(p,\gamma(p))\le R\}$.

First, we show that $\frac \Gamma \lambda(R)\subseteq \frac{\Gamma(R)\Lambda}\Lambda$. For any $\bar \gamma\in \frac \Gamma \Lambda(R)$, $d_{\bar X}(\bar p,\bar \gamma(\bar p))\le R$, let $\gamma\in \Gamma$ such that
$\pi(\gamma)=\bar \gamma$. We shall show that there is $h_\gamma\in \Lambda$ such that
$\gamma h_\gamma\in \Gamma(R)$, thus $\pi(\gamma h_\gamma)=\bar \gamma\in \frac{\Gamma(R)\Lambda}\Lambda$. First, we may assume $h_1, h_2\in \Lambda$ such that $d_{\bar X}(\bar p,\bar \gamma \bar p)=d_X(h_1p,h_2\gamma(p))$. Because $\Lambda$ is a normal subgroup of $\Gamma$,
$$\split R&\ge d_{\bar X}(\bar p,\bar \gamma(\bar p))=d_X(h_1p,h_2\gamma (p))=d_X(p, h_1^{-1}h_2\gamma (p))\\&=d_X(p,\gamma (\gamma^{-1}h^{-1}h_2\gamma)(p))=d_X(p,\gamma h_\gamma (p)).\endsplit$$

Secondly, we show that $\frac{\Gamma(R)\Lambda}\Lambda\subseteq \frac \Gamma \Lambda(R)$. Let $\bar \gamma\in \frac{\Gamma(R)\Lambda}\Lambda$, $\bar \gamma=\pi(\gamma)$, $\gamma\in \Gamma(R)$. Then
$$d_{\bar X}(\bar \gamma (\bar p),\bar p)=d_{\bar X}(\Lambda(\gamma (p),\Lambda(p))=d_X(\gamma (p),\Lambda(p))\le d_X(\gamma (p),p)\le R,$$
thus $\bar \gamma\in \frac \Gamma \Lambda(R)$
\qed\enddemo

\enddocument